\documentclass{amsart}

\usepackage{amsmath,amssymb,url,color,mathtools}
\usepackage[margin=4cm]{geometry}
\usepackage{tikz}
\usepackage{tikz-cd}
\usetikzlibrary{arrows.meta,matrix,positioning,fit}

\newtheorem{lemma}{Lemma}[section]
\newtheorem{prop}[lemma]{Proposition}
\newtheorem{cor}[lemma]{Corollary}
\newtheorem{thm}[lemma]{Theorem}
\newtheorem{example}[lemma]{Example}
\newtheorem{thm?}[lemma]{Theorem?}

\newtheorem{remark}[lemma]{Remark}

\newcommand{\cp}[1]{\tikz[baseline=(char.base)]{\node[shape=circle,draw,inner sep=2pt](char){\smash{\raisebox{-2.5pt}{$\kern .5pt p^{#1}$}}};}
}

\begin{document}
\title[The Group-Theoretic Prime Ax-Katz Theorem]%
{Functional degrees and arithmetic applications III: Beyond Prime Exponent}

\author{Pete L.\ Clark}
\author{Uwe Schauz}

\renewcommand\atop[2]{\genfrac{}{}{0pt}{}{#1}{#2}}
\newcommand{\Mod}[1]{\ (\mathrm{mod}\ #1)}
\newcommand{\etalchar}[1]{$^{#1}$}
\newcommand{\F}{\mathbb{F}}
\newcommand{\et}{\textrm{\'et}}
\newcommand{\ra}{\ensuremath{\rightarrow}}
\newcommand{\lra}{\ensuremath{\longrightarrow}}
\newcommand{\FF}{\F}
\newcommand{\ff}{\mathfrak{f}}
\newcommand{\Z}{\mathbb{Z}}
\newcommand{\N}{\mathbb{N}}
\newcommand{\NN}{\widetilde{\N}}
\newcommand{\mm}{\underline{m}}
\newcommand{\nn}{\underline{n}}
\newcommand{\ch}{}
\newcommand{\R}{\mathbb{R}}
\renewcommand{\P}{\mathbb{P}}
\newcommand{\PP}{\mathbf{P}}
\newcommand{\pp}{\mathfrak{p}}
\newcommand{\C}{\mathbb{C}}
\newcommand{\Q}{\mathbb{Q}}
\newcommand{\ab}{\operatorname{ab}}
\newcommand{\Aut}{\operatorname{Aut}}
\newcommand{\gk}{\mathfrak{g}_K}
\newcommand{\gq}{\mathfrak{g}_{\Q}}
\newcommand{\OQ}{\overline{\Q}}
\newcommand{\Out}{\operatorname{Out}}
\newcommand{\End}{\operatorname{End}}
\newcommand{\Gal}{\operatorname{Gal}}
\newcommand{\CT}{(\mathcal{C},\mathcal{T})}
\newcommand{\lcm}{\operatorname{lcm}}
\newcommand{\Div}{\operatorname{Div}}
\newcommand{\OO}{\mathcal{O}}
\newcommand{\rank}{\operatorname{rank}}
\newcommand{\tors}{\operatorname{tors}}
\newcommand{\IM}{\operatorname{IM}}
\newcommand{\CM}{\mathbf{CM}}
\newcommand{\HS}{\mathbf{HS}}
\newcommand{\Frac}{\operatorname{Frac}}
\newcommand{\Pic}{\operatorname{Pic}}
\newcommand{\coker}{\operatorname{coker}}
\newcommand{\Cl}{\operatorname{Cl}}
\newcommand{\loc}{\operatorname{loc}}
\newcommand{\GL}{\operatorname{GL}}
\newcommand{\PGL}{\operatorname{PGL}}
\newcommand{\PSL}{\operatorname{PSL}}
\newcommand{\Frob}{\operatorname{Frob}}
\newcommand{\Hom}{\operatorname{Hom}}
\newcommand{\Coker}{\operatorname{\coker}}
\newcommand{\Ker}{\ker}
\newcommand{\g}{\mathfrak{g}}
\newcommand{\sep}{\operatorname{sep}}
\newcommand{\new}{\operatorname{new}}
\newcommand{\Ok}{\mathcal{O}_K}
\newcommand{\ord}{\operatorname{ord}}
\newcommand{\Ohell}{\OO_{\ell^{\infty}}}
\newcommand{\cc}{\mathfrak{c}}
\newcommand{\ann}{\operatorname{ann}}
\renewcommand{\tt}{\mathfrak{t}}
\renewcommand{\cc}{\mathfrak{a}}
\renewcommand{\aa}{\mathbf{a}}
\newcommand\leg{\genfrac(){.4pt}{}}
\renewcommand{\gg}{\mathfrak{g}}
\renewcommand{\O}{\mathcal{O}}
\newcommand{\Spec}{\operatorname{Spec}}
\newcommand{\rr}{\mathfrak{r}}
\newcommand{\rad}{\operatorname{rad}}
\newcommand{\SL}{\operatorname{SL}}
\newcommand{\fdeg}{\operatorname{fdeg}}
\renewcommand{\rank}{\operatorname{rank}}
\newcommand{\Int}{\operatorname{Int}}
\newcommand{\zz}{\mathbf{z}}
\newcommand\overbar[1]{\overbracket[.5pt][.7pt]{#1}}

\begin{abstract}
Continuing our work on group-theoretic generalizations of the prime Ax-Katz Theorem,
we give a lower bound on the $p$-adic divisibility of the cardinality of the set of
simultaneous zeros $Z(f_1,f_2,\ldots,f_r)$ of $r$ maps $f_j:A\rightarrow B_j$ between
arbitrary finite commutative groups $A$ and $B_j$ in terms of the invariant factors of
$A, B_1,B_2,\dotsc,B_r$ and the \emph{functional degrees} of the maps
$f_1,f_2,\dotsc,f_r$.
\end{abstract}

\maketitle


\section{Introduction and Main Results}
\noindent
\subsection{Notation and Terminology}\label{sec.term}
Throughout this paper, $p$ is a fixed but arbitrary prime number. We denote by
$\ord_p$ the $p$-\emph{adic valuation} on $\Q$.  For a an integer $q \geq 2$ and a
nonzero integer $Z$, we denote by $\ord_q(Z)$ the largest power of $q$ that divides
$Z$, and we put $\ord_q(0) = \infty$.   (When $q = p$ is prime, this is the $p$-adic
valuation.)  Also we set
$$
\Z^+ :=\, \{n\in\Z\mid n>0\}\ ,\ \
\N\,:=\, \{n\in\Z\mid n\geq0\}\quad\text{and}\quad
\tilde{\N} \,:=\, \N \cup \{ - \infty, \infty\}\,,
$$
and we endow $\tilde{\N}$ with the total ordering that extends the usual ordering on
$\N$ so that $-\infty$ is the smallest and $\infty$ is the largest element.
\\ \\
If $R,R_1, \dotsc, R_r$ are sets and $f_1:R\ra R_1$, \dots, $f_r:R\ra R_r$ are functions
(possibly given as polynomials), such that each of the sets $R_j$ contains a
distinguished element that is denoted as $0=0_{R_j}$, then we define
\[ Z(f_1,\dotsc,f_r)\,=\,Z_{R}(f_1,\dotsc,f_r)
 \,:=\, \bigl\{x \in R \mid f_1(x) = 0\,, \dotsc, f_r(x) = 0\bigr\}. \]
\\
For arbitrary commutative groups $A$ and $B$, we denote with $B^A$ the set of all
functions $f: A \ra B$, and define for each $a \in A$ a \textbf{difference operator}
$\Delta_a \in \End(B^A)$ by
\[\Delta_a f: x \longmapsto f(x+a)-f(x). \]
Following Aichinger-Moosbauer each $f \in B^A$ has a \textbf{functional degree}
$$\fdeg(f)\,:=\,\sup\bigl\{n\in\N\mid \exists a_1 \in A,\dotsc,\exists a_{n} \in A:\,
\Delta_{a_1} \dotsm\Delta_{a_{n}} f\neq0\bigr\}\,\in\,\NN\,,$$
where\footnote{In \cite{Aichinger-Moosbauer21} $\fdeg(0) := 0$\,, and
$\fdeg(f)\in\N\cup\{\infty\}$ for all $f\in B^A$, but we set $\fdeg(0) := -\infty$.}
$\sup(\emptyset):=-\infty$. This degree may be infinite but we fucus mainly on the
subset\footnote{$\mathcal{F}(A,B)$ is actually a $\Z[A]$-submodule of the
\(\Z[A]\)-module $B^A$ over the integral group ring $\Z[A]$ of $A$, where the scalar
multiple of $f\in B^A$ by $c=\sum_{a\in A}n_a[a]\in\Z[A]$ is $cf:x\longmapsto\sum_{a\in
A}n_af(x+a)$.}
\[ \mathcal{F}(A,B) \,:=\, \{ f \in B^A \mid \fdeg(f) < \infty \}.\]

\subsection{The Story so far}
This paper is a direct continuation of our prior works \cite{Clark-Schauz22} and
\cite{Clark-Schauz23a}; in these papers as well as in the present paper, our goal is to
synthesize, further develop and apply work of Wilson \cite{Wilson06} and
Aichinger-Moosbauer \cite{Aichinger-Moosbauer21}.
\\ \\
In \cite{Aichinger-Moosbauer21}, Aichinger and Moosbauer develop a calculus of finite
differences (discrete derivatives) for maps $f: A \ra B$ between arbitrary commutative
groups $A$ and $B$, and introduce the functional degree $\fdeg(f)$ based on the
simple idea that the functional degree should decrease if a discrete derivative is taken.
See also the exposition in our prior work \cite[\S 2.3]{Clark-Schauz22}, where difference
operators and other basics are introduced in more detail.  One of the key insights of
\cite{Aichinger-Moosbauer21} is that it is often fruitful to view the elements of
$\mathcal{F}(A,B)$ as the ``polynomial functions from $A$ to $B$'', a point of view that
was introduced in \cite{Schauz14} already. If $P$ is a polynomial expression in $n$
variables with coefficients in a (not necessarily commutative) rng\footnote{Not a typo: a
ring has a multiplicative identity, a rng may not.} $R$, and $E(P): R^n \!\ra R$ is the
corresponding polynomial function, then, by \cite[Lemma
12.5]{Aichinger-Moosbauer21},
\[ \fdeg(E(P)) \,\leq\, \deg(P)\,. \]
It is an interesting problem to precisely understand the discrepancy between these
two kinds of degree.  After work of Aichinger-Moosbauer \cite[\S
10]{Aichinger-Moosbauer21} and work of the present authors \cite[Prop.\ 2.19 and
Thm.\ 4.9]{Clark-Schauz23a}, we know how to compute $\fdeg(E(P))$ from the family
of coefficients of the monomials of $P$ when $R$ is a commutative integral domain.
Equality holds without restrictions, for all polynomials $P$, if and only if $R$ has
characteristic $0$. In general, however, the functional degree may even be limited by
a constant. For commutative groups $A$ and $B$, we put
\[ \delta(A,B) \,\coloneqq\, \sup_{f \in B^A}  \fdeg(f). \]
When $A$ and $B$ are nontrivial finite commutative groups, Aichinger-Moosbauer
showed that $\delta(A,B)$ is finite if and only if $A$ and $B$ are both $p$-groups for
the same prime number $p$, and they raised the question of determining the exact
value of $\delta(A,B)$ in this case. This was answered by the present authors in
\cite[Thm.\ 4.9\,c]{Clark-Schauz22}, \cite{Schauz21}, and \cite[Thm.\ 3.9]{Schauz14},
using arithmetic results of Weisman \cite{Weisman77} and Wilson \cite{Wilson06}:
\begin{thm}
\label{MAXFDEGTHM} Let $N,\beta,\alpha_1,\ldots,\alpha_N \in \Z^+\!$, let
$\underline{\alpha}:=(\alpha_1,\ldots,\alpha_N )$, and let $B$ be a finite
commutative $p$-group of exponent $p^{\beta}$\!.  Then
\[\delta\bigl( \bigoplus_{i=1}^N \Z/p^{\alpha_i} \Z, B\bigr) \,=\, \delta_p(\underline{\alpha},\beta)
 \]
where
\[\delta_p(\underline{\alpha},\beta)\,\coloneqq\,
 \sum_{i=1}^N( p^{\alpha_i}\!-1)+(\beta-1)(p-1) p^{\max\{\alpha_1,\dotsc,\alpha_N\}-1}.\]
\end{thm}

\noindent For finite commutative $p$-groups $A$ and $B$, the quantity $\delta(A,B)$
can be interpreted as the largest possible ``complexity'' for a map $f: A \ra B$.  For
instance, if $A = (\Z/p\Z)^n$ and $B = \Z/p\Z$, then (as Aichinger-Moosbauer knew) the
largest possible functional degree is $(p-1)n$, and one function of this degree is given
by evaluating the polynomial $t_1^{p-1} \cdots t_n^{p-1}$.  This is related to an
observation of Chevalley: over a finite field $\F_q$, the function $x \mapsto x^q - x$ is
identically zero, so for any polynomial $P \in \F_q[t_1,\ldots,t_n]$ there is another
\textbf{reduced polynomial} $\kern 1pt \overline{\kern-1pt P}\in\F_q[t_1,\ldots,t_n]$
consisting of monomial terms $t_1^{a_1} \cdots t_n^{a_n}$ with $0 \leq a_i \leq q-1$
and such that $E(P) = E(\kern 1pt \overline{\kern-1pt P})$, i.e., the two polynomials
determine the same polynomial function. The largest degree of a reduced monomial is
therefore $(q-1)n$.
\\ \\
Already this hints that the Aichinger-Moosbauer functional calculus should have
numerous theoretic connections, in particular to the following celebrated results.

\begin{thm}\label{CWAK}
Let $q:=p^N$\!.  Let $f_1,\dotsc,f_r\in\F_q[t_1,\dotsc,t_n]$ be polynomials of positive
degrees.  If $\sum_{j=1}^r \deg(f_j)< n$\,, then
\begin{itemize}
\item[a)] $\,\ord_p\bigl(\# Z_{\F_q^n}(f_1,\dotsc,f_r)\bigr) \geq 1$
    \hfill\emph{(Chevalley-Warning Theorem \cite{Chevalley35}, \cite{Warning35})},
\item[b)] $\,\ord_q\bigl(\# Z_{\F_q^n}(f_1,\dotsc,f_r)\bigr) \geq \Bigl{\lceil}
    \frac{n-\sum_{j=1}^r \deg(f_j)}{\max_{j=1}^r \deg(f_j)} \Bigr{\rceil}$
    \hfill\emph{(Ax-Katz Theorem \cite{Ax64}, \cite{Katz71})}.
\end{itemize}
\end{thm}
\noindent Indeed, Aichinger-Moosbauer used their functional calculus to prove the
following result:

\begin{thm}(Group-Theoretic Chevalley-Warning Theorem)
\label{GTCW} Let
\[ A \,\coloneqq\, \bigoplus_{i=1}^m \Z/p^{\alpha_i} \Z\,,\ \ B
 \,\coloneqq\, \bigoplus_{i=1}^n \Z/p^{\beta_i}\Z \]
be finite commutative $p$-groups, and let $f_1,\dotsc,f_r: A^N\!\ra B$ be nonzero
functions. If
\begin{equation*}
\label{GTCWEQ}
\biggl( \sum_{j=1}^r \fdeg(f_j) \biggr) \biggl( \sum_{i=1}^n (p^{\beta_i}\!-1) \biggr)
\,<\, N \sum_{i=1}^m (p^{\alpha_i}\!-1)
\end{equation*}
then
\[ \ord_p(\# Z_{A^N}(f_1,\ldots,f_r)) \,\geq\, 1\,. \]
\end{thm}
\begin{proof}
This is \cite[Thm.\ 12.2]{Aichinger-Moosbauer21} with zero functions of degree
$-\infty$ excluded.
\end{proof}
\noindent
If $R$ is a finite rng of prime power order and $P_1,\ldots,P_r$ are polynomial
expressions over $R$ in $N$ variables, then applying Theorem \ref{GTCW} with $A = B
= (R,+)$ and with $f_1 = E(P_1),\ldots,f_r = E(P_r)$ the associated polynomial
functions from $R^N$ to $R$, one obtains a ring-theoretic generalization of of Theorem
\ref{CWAK}\,a.  This gives a ``psychological'' explanation for the presence of $N$ in
Theorem \ref{GTCW}, but nothing is lost by taking $N = 1$.
\\ \\
The same work \cite{Aichinger-Moosbauer21} gave a group-theoretic generalization of
Warning's Second Theorem \cite[Thm.\ 14.2]{Aichinger-Moosbauer21},  but they left
open the problem of applying their calculus to higher $p$-adic congruences. However, a
2006 work of R.\ Wilson \cite{Wilson06} gave a strikingly new and elementary proof of
Theorem \ref{CWAK}\,b over the prime field $\F_p$ using, in particular, the difference
operators $\Delta: f \mapsto\Delta f:=\bigl(x\mapsto f(x+1)-f(x)\bigr)$ from the calculus
of finite differences. Comparing the work of Wilson with that of Aichinger-Moosbauer,
we found that -- notwithstanding some differences in perspective and presentation --
they are deeply related.  Our proof of Theorem \ref{MAXFDEGTHM} makes use either
of Wilson's work or, alternately, earlier related work of Weisman \cite{Weisman77}.
Moreover, with some further development of the Achinger-Moosbauer calculus --
especially that for commutative group $B$, the elements of $\mathcal{F}(\Z^N\!,B)$
have series expansions (as recalled in Theorem \ref{4.2}) -- we were able \cite[Cor.\
1.9]{Clark-Schauz23a} to refine Wilson's argument to give the following group-theoretic
generalization of Theorem \ref{CWAK}\,b over $\F_p$:

\begin{thm}[Group-Theoretic Prime Ax-Katz Theorem]
\label{CS2MAIN1} Let $N,n,r \in \Z^+$\!, and put $A := (\Z/p\Z)^N$\!. Let
$f_1,\dotsc,f_r \in A^{A^n}$ be nonconstant functions. Then
\[\ord_p(\# Z_{A^n}(f_1,\dotsc,f_r)) \,\geq\,
\biggl{\lceil} \frac{N\bigl(n-\sum_{j=1}^r \fdeg(f_j)\bigr)}
{\max_{j=1}^r \fdeg(f_j)}\biggr{\rceil}.\]
\end{thm}
\noindent We emphasize that like Theorem \ref{GTCW}\,, Theorem \ref{CS2MAIN1}
is a purely group-theoretic result.  When $A = \F_p$ we recover Theorem
\ref{CWAK}\,b over the prime field $\F_p$.  When $A = \F_q$, because of the
connection between the functional degree and the $p$-weight degree, it recovers
Moreno-Moreno's strengthening of the prime Ax-Katz Theorem
\cite{Moreno-Moreno95}, which however does not imply the full Ax-Katz Theorem
over $\F_q$ (cf.\ \cite[Remark 1.4]{Clark-Schauz23a}).
\\ \\
After seeing a related manuscript of Grynkiewicz \cite{Grynkiewicz22}, we noticed
that the argument that proves Theorem \ref{CS2MAIN1} can be adapted to prove a
more general result:


\begin{thm}\cite[Thm.\ 1.7]{Clark-Schauz23a}
\label{CS2MAIN2} Let $N,r,\beta_1,\dotsc,\beta_r  \in \Z^+$\!, and put $A :=
(\Z/p\Z)^N$\!. For each $1 \leq j \leq r$, let $f_j \in (\Z/p^{\beta_j} \Z)^A$ be a nonzero
function. Then
\[ \ord_p(\# Z_{A}(f_1,\dotsc,f_r)) \,\geq\, \biggl\lceil \frac{N- \sum_{j=1}^r \frac{p^{\beta_j}-1}{p-1} \fdeg(f_j)}
 {\max_{j=1}^r \,p^{\beta_j-1} \fdeg(f_j)} \biggr\rceil. \]
\end{thm}

\noindent Let us compare Theorems \ref{CS2MAIN1} and \ref{CS2MAIN2}.  In the
former result all the maps $f_1,\ldots,f_r$ take values in a fixed finite commutative
$p$-group $A$ that is required to have exponent $p$ but is not required to be cyclic.
In the latter result the maps $f_1,\ldots,f_r$ take values in varying cyclic $p$-groups
$\Z/p^{\beta_j} \Z$ that are not necessarily of prime exponent.  So it may seem that
we have lost generality in the requirement that the target groups be cyclic.  But there
is a \textbf{cyclic exchange} phenomenon: as discussed in Section \ref{S2.1},  in
these results we may exchange a map $f_j$ into a finite commutative $p$-group with
$K$ invariant factors for a $K$-tuple $(f_{j,1},\ldots,f_{j,K})$ of maps into finite cyclic
$p$-groups compatibly with our setup.  Thus considering maps with targets in
\emph{varying cyclic groups} carries all the content of the general case.
\\ \\
Moreover, there is a Sylow primary decomposition for maps of finite functional degree
between arbitrary finite commutative groups \cite[Cor.\ 3.14\,c \& Cor.\
3.15]{Clark-Schauz22}.  Using this, one can extend all of these group-theoretic results
from finite commutative $p$-groups to arbitrary finite commutative groups. This was
done already in \cite[Rem.\ 1.8 \& Cor.\ 1.9]{Clark-Schauz23a},  and we repeat the
discussion here in Section \ref{S2.2}.





\subsection{The Main Theorem}
The main result of this paper is an Ax-Katz type lower bound on $\ord_p(\#
Z(f_1,\ldots,f_r)$ for maps $f_j: A \ra B_j$ between arbitrary finite commutative
$p$-groups $A,B_1,\ldots,B_r$.   As above, by ``cyclic exchange'' (cf.\ Section
\ref{S2.1}) we may assume that each $B_j$ is cyclic, so we may write
$$A\,=\,\bigoplus_{i=1}^N \Z/p^{\alpha_i} \Z\,,$$
and
$$B_1=\Z/p^{\beta_1}\Z\,,\,\dotsc\,,\, B_r=\Z/p^{\beta_r}\Z\,,$$
where $r,\beta_1,\dots,\beta_r ,N,\alpha_1,\dotsc,\alpha_N\in \Z^+\!$. We may assume
without loss of generality that each $f_j$ is nonconstant of functional degree at most
$d_j\in\Z^+\!$. In other words, for each $1\leq j\leq r$ we have a function
$$f_j: A \ra B_j\quad\text{with}\quad0\,<\,\fdeg(f_j)\,\leq\,d_j\,.$$
Put
$$\mathcal{A}\,:=\,\sum_{i=1}^N\frac{p^{\alpha_i}\!-1}{p-1}
\quad\text{and}\quad
\mathcal B\,:=\,\sum_{j=1}^rd_j\frac{p^{\beta_j}\!-1}{p-1}\,.$$
We may order the $\alpha_i$ and the $f_j$ so that
$$\quad\alpha_1 \geq \alpha_2 \geq \dotsb \geq\alpha_N\, \ \quad\text{and }
d_1p^{\beta_1}\geq d_2p^{\beta_2}\geq\dotsb\geq d_rp^{\beta_r}.\smallskip$$

\noindent To express our result we also need the \textbf{conjugates}
$\alpha_1',\alpha_2',\dotsc,\alpha_{\alpha_1}'$ defined by
$$\alpha_j'\,:=\,\#\bigl\{1\leq i\leq N\mid\alpha_i\geq j\bigr\}\,,$$
which we discuss in more detail in Section \ref{sec.CS}. We set
\begin{align*}
\alpha\,:={}&\,\alpha_1+ \alpha_2+ \dotsb + \alpha_N\\
 \,={}&\,\alpha_1'+ \alpha_2'+ \dotsb + \alpha_{\alpha_1}'\quad\text{(by Example \ref{Ex.con1})}\\
\intertext{and}
\breve\alpha\,:={}&\,\breve\alpha_1+ \breve\alpha_2+ \dotsb + \breve\alpha_N\\
 \,={}&\,\alpha_1'+ \alpha_2'+ \dotsb + \alpha_{\breve\alpha_1}'\quad\text{(by Example \ref{Ex.con2})}\\
\intertext{where}
\breve\alpha_i\,:={}&\,\min\bigl\{\alpha_i\,,L\bigr\}\\
\intertext{with}
L\,:={}&\,\beta_1+\lfloor\log_p\bigl(d_1\bigr)\rfloor\,.\\
\intertext{Using that $\alpha=\alpha_1'+ \alpha_2'+ \dotsb + \alpha_{\alpha_1}'$\,, we
define numbers $D_1,D_2,\dotsc,D_{\alpha}$ by setting}
\bigl(D_1,D_2,\dotsc,D_{\alpha}\bigr)
\,:={}&\,\bigl(\,\underbrace{1,1,\dotsc,1}_{\alpha_1'\ \text{times}},
\,\underbrace{p,p,\dotsc,p}_{\alpha_2'\ \text{times}},
\,\dotsc,\,\underbrace{p^{\alpha_1-1},p^{\alpha_1-1},\dotsc,p^{\alpha_1-1}}%
 _{\alpha_{\alpha_1}'\ \text{times}}\,\bigr).\\
\intertext{We further put}
\breve{\mathcal{A}}\,:={}&\,\sum_{i=1}^N\frac{p^{\breve\alpha_i}\!-1}{p-1}\\
\,={}&\,\alpha_1'p^0+\dotsb+\alpha_{\breve\alpha_1}'p^{\breve\alpha_1-1}
\quad\text{(by Example \ref{Ex.con2})}\\[5pt]
 \,={}&\,D_1+\dotsb+D_{\breve\alpha}
\quad\text{(as $\breve\alpha=\alpha_1'+\alpha_2'+\dotsb+\alpha_{\breve\alpha_1}'$
by Example \ref{Ex.con2}).}\\
\end{align*}

\noindent With these parameters and definitions we can prove in the subsequent
sections and in particular in Section \ref{sec.proof}, the following main result of our
paper, which is a simultaneous generalization of Theorem \ref{GTCW} of
Aichinger-Moosbauer and Theorem \ref{CS2MAIN1} (hence also of Theorem
\ref{CS2MAIN2}) of the present authors.

\begin{thm}
\label{MAINTHM}
With the parameters and settings above,
$$
\ord_p \bigl(\# Z_{A}(f_1,\dotsc,f_r))
 \,\geq\,\begin{cases}
 \displaystyle{\biggl\lceil\frac{\breve{\mathcal{A}}-\mathcal B}{d_1p^{{\beta_1-1}}}\biggr\rceil
 +\alpha-\breve\alpha}
 & \text{if $\breve{\mathcal{A}}>\mathcal B$,}\\
 \,\alpha-\max\bigl\{1\leq t\leq\alpha\mid D_1+\dotsb+D_{t}\leq \mathcal B\bigr\}^{\strut}\!
 & \text{if $\breve{\mathcal{A}}\leq\mathcal B$.}
\end{cases}
$$\smallskip
\end{thm}

\noindent Note that our lower bound is equal to $0$ if $\mathcal{A}\leq\mathcal B$.
This is because $\breve{\mathcal{A}}\leq\mathcal{A}\leq\mathcal{B}$ (so we are in
the second case of Theorem \ref{MAINTHM}), and by Example \ref{Ex.con1} also
$$D_1+\dotsb+D_{\alpha}
 \,=\,\alpha_1'p^0+\dotsb+\alpha_{\alpha_1}'p^{\alpha_1-1}
 \,=\,\sum_{i=1}^N\frac{p^{\alpha_i}\!-1}{p-1}
 \,=\,\mathcal{A}\,\leq\,\mathcal{B}\,.$$
In the case $\mathcal{A}>\mathcal B$, however,
$D_1+\dotsb+D_{\alpha}=\mathcal{A}\nleq\mathcal{B}$ and the lower bound in
Theorem \ref{MAINTHM} is positive (in both cases), so that we obtain the following
corollary:

\begin{cor}\label{MAINCOR1}
If $\mathcal{A}>\mathcal B$ then
\[ \ord_p\bigl(\# Z_{A}(f_1,\ldots,f_r)\bigr) \,\geq\, 1.\smallskip \]
\end{cor}
\noindent As we explain in Section \ref{S2.1}, Corollary \ref{MAINCOR1} implies
Theorem \ref{GTCW}: essentially, it is the generalization of  Theorem \ref{GTCW} in
which the maps $f_1,\ldots,f_r$ are allowed to take values in varying finite commutative
$p$-groups $B_1,\ldots,B_r$.\footnote{Aichinger-Moosbauer's proof of Theorem
\ref{GTCW} can be adapted to prove this generalization.}
\\ \\
\noindent Our main theorem takes a somewhat simpler form when $\alpha_1 = \dotsb =
\alpha_N$. If we set
$$Q\,:=\,\bigl\lfloor\log_p\bigl((p-1)\mathcal B/N+1\bigr)\bigr\rfloor\quad\text{and}\quad
 R\,:=\,\biggl\lfloor\frac{\mathcal B-N\frac{p^Q-1}{p-1}}{p^Q}\biggr\rfloor$$
then, as we prove in Section \ref{sec.proof2}, we obtain the following corollary:

\begin{cor}\label{MAINCOR}
If $\alpha_1=\dotsb=\alpha_N$ then
$$
\ord_p \bigl(\# Z_{A}(f_1,\dotsc,f_r))
 \,\geq\,\begin{cases}
 \displaystyle{\biggl\lceil\frac{N\frac{p^{\breve\alpha_1}-1}{p-1}-\mathcal B}{d_1p^{{\beta_1-1}}}\biggr\rceil}
 \,+\,N(\alpha_1-\breve\alpha_1)
 & \text{if $N\frac{p^{\breve\alpha_1}-1}{p-1}>\mathcal B$,}\\
 \,N(\alpha_1-Q)-R^{\strut\phantom{|}}
 & \text{if $N\frac{p^{\breve\alpha_1}-1}{p-1}\leq\mathcal B$.}
\end{cases}
$$
\end{cor}
\noindent Theorem \ref{CS2MAIN2} is the case $\alpha_1 = \ldots = \alpha_N = 1$ of
Corollary \ref{MAINCOR} (with $d_j=\fdeg(f_j)$ for $1\leq j\leq r$). Indeed, in this case
$\breve\alpha_1=1$.  If  $N > \mathcal{B}$ then the lower bounds in Theorem
\ref{CS2MAIN2} and Corollary \ref{MAINCOR} coincide, whereas if $N \leq
\mathcal{B}$ neither lower bound is positive so the results are vacuous in this case.

\subsection{Schedule of Remaining Tasks}

In Section 2, we explain how our results can be applied in more general situations and
to polynomials over rngs. First in Section 2.1 we look at groups of prime power order,
not necessarily cyclic ones. In Section 2.2, we generalize then to commutative groups of
finite order. In Section 2.3, polynomials over rngs of prime power order are discussed.
In Section 2.4, this is generalized to polynomials over rngs of finite order.

In Sections 3.4 and 3.5, we present the proof of Theorem \ref{MAINTHM} and Corollary
\ref{MAINCOR}. The proof of Theorem \ref{MAINTHM} actually requires quite some
preparation, and the more obvious part of that preparation is given in Sections 3.1
through 3.3.

There are also two less obvious optimization tasks. The necessity to study these kind
of optimizations becomes apparent only during the main proof, after the functions
$\nu_p(\protect\underline{\alpha},\bullet)$ and $\mathcal N$ where introduced in
Section 3.2 and 3.4. In other words, the main proof motivates and sets up those
tasks. We moved those two optimization tasks into the subsequent Sections 4 and 5,
as they can be studied independently. Inside our main proof in Section 3.4, the
results of those investigations are then just cited.

Section 4 is about the second lower bound of Theorem \ref{MAINTHM}:
$$\alpha-\max\bigl\{1\leq t\leq\alpha\mid D_1+\dotsb+D_{t}\leq \mathcal B\bigr\}\,.$$
It includes in Section 4.2 a supplementary discussion of alternative ways to express
this lower bound, which is not needed to understand the main proof and the rest of
the paper.

Section 5 deals with the case distinction, and the first lower bound of Theorem
\ref{MAINTHM}: $$\biggl\lceil\frac{\breve{\mathcal{A}}-\mathcal
B}{d_1p^{{\beta_1-1}}}\biggr\rceil+\alpha-\breve\alpha,.$$

Section 6 provides some background on conjugate partitions, with a number of
lemmas and examples that we conveniently cite throughout the
paper. 
Readers not familiar with conjugate partitions and Ferrers diagrams may want to read
through this section
first.

\section{Extending the Scope: Finite Commutative Groups and Polynomials}\label{S2}

In this section we explain how our results can be applied in more general
situations. 

\subsection{Commutative Groups of Prime Power Order}\label{S2.1}

As mentioned before, our results can still be applied when arbitrary finite commutative
$p$-groups $B_j=\bigoplus_{i=1}^{K_j} (\Z/p^{\beta_{j,i}} \Z)$ replace the cyclic
$p$-groups $\Z/p^{\beta_j} \Z$ as codomains. One just has to use the coordinate
projections $\pi_k: \bigoplus_{i=1}^{K_j} (\Z/p^{\beta_{j,i}}\Z) \ra \Z/p^{\beta_{j,k}} \Z$
first (as already explained in \cite[Rem.\ 1.8]{Clark-Schauz23a}), to define the functions
$$f_{j,k}\,:=\,\pi_k \circ f_j\quad\text{with}\quad\fdeg(f_{j,k})\leq \fdeg(f_j)\leq d_j\,.$$
Then our results can be applied to those $f_{j,k}$ and the degree restrictions
$\fdeg(f_{j,k})\leq d_j$\,, with the outcome interpreted in terms of the $f_j$.

If, for instance, Corollary \ref{MAINCOR1} is applied in this way, then the parameter
$\mathcal  B$ takes the form $\mathcal B=\sum_{j=1}^r
\bigl(d_j\sum_{k=1}^{K_j}\frac{p^{\beta_{j,k}}\!-1}{p-1}\bigr)$, while the parameter
$\mathcal A$ and the conclusion remain unchanged: still $\mathcal A>\mathcal B$
implies $\ord_p(\# Z_{A}(f_1,\ldots,f_r)) \geq 1$. This generalizes Theorem \ref{GTCW},
where $B_1=\dotsb=B_r=\bigoplus_{i=1}^{K} (\Z/p^{\beta_{i}} \Z)$, i.e., $\mathcal
B=\sum_{j=1}^r d_j\sum_{k=1}^{K}\frac{p^{\beta_{k}}\!-1}{p-1}$, and
$A=\bigl(\bigoplus_{i=1}^m \Z/p^{\alpha_i}\Z\bigr)^N$\!\!, i.e.,
$\mathcal{A}=N\sum_{i=1}^m\frac{p^{\alpha_i}\!-1}{p-1}$.\smallskip

\subsection{Commutative Groups of Finite Order}\label{S2.2}

Let $A,B_1,\ldots,B_r$ be any nontrivial finite commutative groups. We write out the
primes dividing $\# \left( A \times \prod_{i=1}^r B_i \right)$ as $\ell_1 < \ldots <
\ell_s$\,, and set
 $$A[\ell_j^{\infty}]:=\{x \in A \mid \ell_j^k x = 0 \text{ for some } k\in\Z^+\}\,.$$
For each fixed $1 \leq j \leq r$, we have a canonical $\Z$-module injection
 \[ \prod_{h=1}^s B_j[\ell_h^{\infty}]^{A[\ell_h^{\infty}]} \ra B_j^{\,A}\]
in which we send each vector $(g_{j,1},\ldots,g_{j,s})$ of functions $g_{j,h}:
A[\ell_h^{\infty}] \ra B_j[\ell_h^{\infty}]$ to the identically named function
\begin{align*}
    (g_{j,1},\ldots,g_{j,s}):A=\prod_{h=1}^sA[\ell_h^{\infty}] &\ra
    B_j=\prod_{h=1}^sB_j[\ell_h^{\infty}]\,,\\
    (x_1,\ldots,x_s) &\mapsto \bigl(g_{j,1}(x_1),\ldots,,g_{j,s}(x_s)\bigr).
\end{align*}
By \cite[Cor.\ 3.14\,c \& Cor.\ 3.15]{Clark-Schauz22}, upon restriction to functions of
finite functional degree, this yields the canonical isomorphy
\[\mathcal{F}(A,B_j) \,=\, \prod_{h=1}^s \mathcal{F}(A[\ell_h^{\infty}],B_j[\ell_h^{\infty}])\]
in which moreover $\fdeg((g_{j,1},\ldots,g_{j,s})) = \max \{\fdeg(g_{j,h})\mid 1\leq h\leq
s\}$. In other words, a map $f_j: A \ra B_j$ of finite functional degree is determined by
its restrictions $g_{j,h}:=f_j|_{A[\ell_h^\infty]}\in B_j[\ell_h^{\infty}]^{A[\ell_h^{\infty}]}$.
We have
$$f_j=(g_{j,1},\ldots,g_{j,s})\quad\text{and}\quad\fdeg(f_j)
 = \max \{\fdeg(g_{j,h})\mid 1\leq h\leq s\}.$$
If we now consider $r$ maps $f_j: A \ra B_j$, where $1 \leq j \leq r$, then we get
$s\times r$ primary component maps $g_{j,h}: A[\ell_h^{\infty}]\ra B_j[\ell_h^{\infty}]$,
and it is immediate that
\[ \# Z(f_1,\ldots,f_r) \,=\, \prod_{h=1}^s \# Z(g_{1,h}\,,\ldots,g_{r,h})\,. \]
So, with the previous Section\,\ref{S2.1}, we obtain for each $1 \leq h \leq s$ a lower
bound on $\ord_{\ell_h}\bigl(\# Z(f_1,\ldots,f_r)\bigr)$ in terms of $A,B_1,\ldots,B_r$
and $\fdeg(f_1),\ldots,\fdeg(f_r)$.

\subsection{Polynomials over Rngs of Prime Power Order}\label{S2.3}

Before \cite{Clark-Schauz23a}, Ax-Katz type $p$-adic congruences on the solution set
of a polynomial system over a finite rng were only known for finite commutative rings in
which every ideal is principal \cite{Ax64}, \cite{Katz71}, \cite{Marshall-Ramage75},
\cite{Katz12}. Now let $R$ be a finite rng with order a power of $p$, so there are
$N,\alpha_1,\ldots,\alpha_N\in\Z^+$ such that
\[ (R,+) \,\cong\, \bigoplus_{i=1}^N\,\Z/p^{\alpha_i} \Z\,=:\,A_1\,. \]
Let $P_1,\ldots,P_r$ be polynomials in $n$ variables over $R$ with $\deg(P_j)\leq d_j$
for each $1\leq j\leq r$\,. Then Theorem \ref{MAINTHM} with $A:=A_1^n$ and the
previous Section \ref{S2.1} apply to give an Ax-Katz type lower bound on $\ord_p(\#
Z(P_1,\ldots,P_r))$. In particular, as $\breve{\mathcal{A}}\geq n$\,, one sees the
following \textbf{asymptotic Ax-Katz over a finite rng}: if $r$ and $d_1,\ldots,d_r$ remain
fixed, then $\ord_p(\# Z(P_1,\ldots,P_r))$ approaches infinity with $n$\,.
\\ \\
This asymptotic Ax-Katz result is also established in a concurrent work by the first
author and N.\ Triantafillou \cite[Thm.\ 6.2]{Clark-Triantafillou23}.  The proof given
there uses a new invariant:  for nontrivial, finite commutative $p$-groups $A$ and $B$
the \textbf{summation invariant} $\sigma(A,B)$ is the largest $d \in \tilde{\N}$ such
that $\int_A f = 0$ for all maps $f: A \ra B$ with functional degree at most $d$.  The
invariant $\sigma(A,B)$ does not appear explicitly in the work of
Aichinger-Moosbauer, but neverthless they give what amounts to a computation of
$\sigma(A,B)$ when $A$ and $B$ are $p$-groups and $B$ has exponent $p$
\cite[Lemma 12.1]{Aichinger-Moosbauer21}, and this is a key ingredient of their proof
of  Theorem \ref{GTCW}.  In \cite{Clark-Triantafillou23} lower bounds are given on
$\sigma(A,B)$ in the general case (and exact computations are given in some further
special cases).  If one takes these results as a ``black box,'' then the proof of
asymptotic Ax-Katz over a finite rng given in \cite{Clark-Triantafillou23} is much
simpler than the proof of our main result.  However the lower bound on $\ord_p(\#
Z(P_1,\ldots,P_r))$ given by our Theorem \ref{MAINTHM} is better than (or equal to,
in certain special cases) the corresponding bound given by the methods of
\cite{Clark-Triantafillou23}. Thus neither work majorizes the other.  \\ \indent Our
present approach implicitly uses a ``lifted variant'' of $\sigma(A,B)$. A comparison of
this lifted variant with $\sigma(A,B)$ is made in \cite[\S 7]{Clark-Triantafillou23}.

\subsection{Polynomials over Rngs of Finite Order}\label{S2.4}

The two previous sections
can be combined to address the case of polynomial expressions in $n$ variables of
degrees $d_1,\ldots,d_r$ over any nontrivial finite rng $R$. In this case, the
asymptotic Ax-Katz phenomenon can be expressed as follows: keeping the number
and degrees of the polynomial expressions $f_1,\ldots,f_r$ fixed, we find that
$\ord_{\# R} (Z(f_1,\ldots,f_r))$ approaches infinity with $n$.   A slightly different, but
equivalent, formulation is given in \cite[Thm.\ 6.2]{Clark-Triantafillou23}.

\section{Reduction to Discrete Optimization}

In this section we prove Theorem \ref{MAINTHM} (and afterwards Corollary
\ref{MAINCOR}) based on some technical results about the minimum values of certain
discrete functions. In other words, this main part of the proof reduces us to some
discrete optimization problems. These discrete optimization problems are stated and
solved completely independent from the original problem in later section, but are cited
and used here. We first recall some basics from our earlier work, then introduce some
basic number theoretic results, and then start that reductionistic proof.

\subsection{Some Recalled Results}
We provide some basics about series expansions in terms of binomial polynomials
$\binom{t}{d} \coloneqq \frac{t(t-1)\dotsm(t-d+1)}{d!} \in \Q[t]$, which are integer valued,
i.e., $\binom{x}{d}\in\Z$ whenever $x\in\Z$. Again, we write $\underline{n}$ for
$(n_1,\dotsc,n_N)$, and we set $|\underline{n}|:=n_1+\dotsb+n_N$.

\begin{thm}
\label{4.2} Let $B$ be a commutative group, and let $f \in B^{\Z^N}$\!\!.
\begin{itemize}
\item[a)] There is a unique function $c: \N^N\!\ra B$ such that
\begin{equation*}
\label{4.2EQ1}
f(\underline{x}) = \sum_{\nn  \in \N^N} \binom{x_1}{n_1} \dotsm
\binom{x_N}{n_N} c(\nn) \quad\text{for all $\underline{x} \in \N^N$\!.}
\end{equation*}
The function values of $c$ are given by the formula $c(\nn) =
\Delta^{\nn}f(\underline{0})$.\smallskip
\item[b)] If $d:=\fdeg(f)< \infty$ then
\begin{equation*}
f(\underline{x}) \,= \sum_{\atop{\nn \in \N^N}{\!\!|\nn| \leq d}} \binom{x_1}{n_1}
\dotsm \binom{x_N}{n_N} \Delta^{\nn}f(\underline{0})\quad\text{for all
$\underline{x} \in \Z^N$\!.}
\end{equation*}
\end{itemize}
\end{thm}
\begin{proof}
This is \cite[Thm.\ 2.8]{Clark-Schauz23a}, where it is also mentioned that the sum in the
first part is always well-defined, because at each fixed point $\underline{x}$ in $\N^N$
(unlike in $\Z^N\setminus\N^N$) the sum has always only a finite number of nonzero
summands, for whatever $c$\,.
\end{proof}

\noindent We now recall some terminology and results concerning proper lifts. Let $\mu:
B \ra B'$ be a surjective homomorphism of commutative groups, and let $f \in
\mathcal{F}(\Z^N\!,B')$. To define a \textbf{proper lift} $\tilde f \in \mathcal{F}(\Z^N\!,B)$
of $f$ (it is not unique), observe first that there exists (by Theorem \ref{4.2}\,b) a unique
(by Theorem \ref{4.2}\,a) function $c: \N^N\!\ra B'$ such that
\[f(\underline{x}) \,= \sum_{\nn  \in \N^N} \binom{x_1}{n_1} \dotsm
 \binom{x_N}{n_N} c(\nn)\quad\text{for all $\underline{x} \in \Z^N$\!.}\]
Then a \textbf{proper lift} of $c$ to $B$ is a function $\tilde c: \N^N\!\ra B$ such that
$$\mu \circ \tilde c \,=\,
 c\quad\text{and}\quad\tilde c(\nn) = 0 \,\iff\, c(\nn) = 0\quad\text{for all $\nn \in \N^N$\!.}$$
Such a proper lift always exists (non-uniquely), and we can use it to define a
\textbf{proper lift} of $f$ by
\[\tilde{f}(\underline{x}) \,:=\, \sum_{\nn  \in \N^N} \binom{x_1}{n_1} \dotsm
 \binom{x_N}{n_N} \tilde c(\nn)\,. \]
This sum is actually well-defined, because $\tilde c: \N^N\!\ra B$ is finitely nonzero (i.e.,
its support is finite), since our $c: \N^N\!\ra B'$ is finitely nonzero (as in Theorem
\ref{4.2}\,b). For every proper lift $\tilde f$ of $f$ we have
$$f\,=\,\mu\circ\tilde{f}\quad\text{and}\quad\fdeg(\tilde{f})\,=\,\fdeg(f)\,.$$
\\
Inside our main proof we work with series expansions of proper lifts of pullbacks, and
we need the following corollary to Theorem \ref{4.2}, in which the maximal finite
degrees
$$\delta_p(\underline{\alpha},h)\,:=\,\sum_{i=1}^N\left( p^{\alpha_i}-1 \right) +
 (h-1)(p-1) p^{\max \{\alpha_1,\dotsc,\alpha_N\} - 1}$$
of Theorem \ref{MAXFDEGTHM} play an important role:

\begin{cor}
\label{COR4.57}  Let $N,\beta,\alpha_1,\dotsc,\alpha_N\in \Z^+$\!. Let $f:
\bigoplus_{i=1}^N \Z/p^{\alpha_i} \Z \ra \Z/p^{\beta}\Z$ be any function, $F: \Z^N\! \ra
\Z/p^{\beta} \Z$ be the pullback of $f$, and $\tilde{F}: \Z^N\! \ra \Z$ be a proper lift of
$F$.
\begin{itemize}
\item[a)]\ \vspace{-.5em}
\begin{equation*}
\tilde F(\underline{x})\
=\!\!\! \sum_{\atop{\nn \in \N^N}{\ \ |\nn| \leq \delta_p(\underline{\alpha},\beta)}}\!\!\!\!
\binom{x_1}{n_1} \dotsm \binom{x_N}{n_N} \Delta^{\nn} \tilde F(\underline{0})
\quad\text{for all $\underline{x} \in \Z^N$\!.}
\end{equation*}
\item[b)] For all $h \in \Z^+\!$ and all $\underline{n} \in \N^N$ with
    $|\underline{n}| > \delta_p(\underline{\alpha},h)$,
    \[ p^h \bigm{|} \Delta^{\nn}\tilde F(\underline{0})\,. \]
\end{itemize}
\end{cor}
\begin{proof} This is \cite[Cor.\ 2.25]{Clark-Schauz23a}.
\end{proof}

\subsection{The numbers \protect{$\nu_p(\protect\underline{\alpha},\protect\underline{n})$}
and \protect{$\mathcal{V}_p(\protect\underline{\alpha},D)$}} For $\alpha \in \Z^+$ and
$n \in \N$\,, we put
\[\nu_p(\alpha,n)\,:=\,\ord_p\Biggl(\,\sum_{x=0}^{p^\alpha\!-1} \binom{x}{n}\Biggr). \]

\begin{lemma}
\label{Lem.vp} For each $\alpha\in \Z^+$\! and $n\in \N$,
$$\nu_p(\alpha,n)
\,=\,\begin{cases}
\alpha-\ord_p(n+1) & \text{if $n\leq p^\alpha\!-1$,}\\
\infty & \text{otherwise.}
\end{cases}$$
\end{lemma}
\begin{proof}
The case $n = 0$ is handled by Proposition \ref{Pr.Vv}\,a, while if $n \geq p^{\alpha}$
then $\sum_{x=0}^{p^\alpha\!-1}\binom{x}{n} = 0$\,, so $\nu_p(\alpha,n) = \infty$. So we
may assume that $1 \leq n \leq p^{\alpha}\!-1$\,.  Using Pascal's rule $\binom{a}{b} =
\binom{a-1}{b} + \binom{a-1}{b-1}$ we see that
\begin{align*}
\sum_{x=0}^{p^\alpha\!-1}\binom{x}{n}
\,\,&=\,\,\binom{n+1}{\,n+1\,}\ +\binom{n+1}{n}+\binom{n+2}{n}+\binom{n+3}{n}+\dotsb+\binom{p^\alpha\!-1}{n}\\
&=\,\,\binom{n+2}{\,n+1\,}\ +\binom{n+2}{n}+\binom{n+3}{n}+\dotsb+\binom{p^\alpha\!-1}{n}\\
&\ \,\vdots\\
&=\,\,\binom{p^\alpha\!-1}{\,n+1\,}\ +\binom{p^\alpha\!-1}{n}\\
&=\,\,\binom{p^\alpha}{\,n+1\,}.
\end{align*}
Now we can apply Kummer's insight \cite{Kummer52} that the $p$-adic valuation of a
binomial coefficient $\binom{a}{b}$ is the number of carries when $b$ and $a-b$ are
added in base $p$\,. Since the base $p$ representation of $p^{\alpha}$ is
$1000\dotsm0$ with $\alpha$ zeros, this number of carries is $\alpha - \ord_p(n+1)$
in our case.
\end{proof}

\noindent For $N\in\Z^+\!$, $\underline{\alpha} = (\alpha_1,\ldots,\alpha_N) \in
(\Z^+)^N$\! and $\underline{n}=(n_1,\ldots,n_N) \in \N^N$\!, we put
\[\nu_p(\underline{\alpha},\underline{n})
 \,:=\,\ord_p\Biggl(\,\sum_{\underline{x}\in[p^{\underline{\alpha}})\!\!}
 \binom{x_1}{n_1}\dotsm\binom{x_N}{n_N}\Biggr)
 \,=\,\ord_p\Biggl(\,\prod_{i=1}^N\sum_{x_i\in[p^{\alpha_i})\!\!}
 \binom{x_i}{n_i}\Biggr)
 =\,\sum_{i=1}^{\smash{N}}\nu_p(\alpha_i,n_i)\,, \]
where
$$[p^{\underline{\alpha}})
 \,:=\,[p^{\alpha_1})\times[p^{\alpha_2})\times\dotsm\times[p^{\alpha_N})
 \quad\text{with}\quad[p^{\alpha_i}):=\{0,1,\dotsc,p^{\alpha_i}\!-1\}\,.$$
To any $D\in\N\cup\{\infty\}$, we also define
$$\mathcal{V}_p(\underline{\alpha},D)
 \,:=\,\min\bigl\{\nu_p(\underline{\alpha},\underline{n})\bigm{|}|\underline{n}|\leq D\bigr\}\,,$$
which is always finite and zero if $D=\infty$, as we see next:

\begin{prop}
\label{Pr.Vv} Let $\alpha \in \Z^+\!$, $\underline{\alpha} \in (\Z^+)^N$\!, and
$D\in\N\cup\{\infty\}$.
\begin{itemize}
\item[a)]\ \vspace{-.5em}
\[ v_p(\alpha,0) \,=\, \alpha\,, \]
and thus
\[ v_p(\underline{\alpha},\underline{0}) \,=\, |\underline{\alpha}|\]
and
\[ \mathcal{V}_p(\underline{\alpha},D) \,\leq\, |\underline{\alpha}|\,.\medskip \]
\item[b)]\ \vspace{-.5em}
\[ v_p(\alpha,p^{\alpha}\!-1) \,=\, 0\,, \]
and thus
\[ v_p\bigl(\underline{\alpha},(p^{\alpha_1}\!-1,\ldots,p^{\alpha_N}\!-1)\bigr) \,=\, 0\]
and
\[ D \geq \sum_{i=1}^N (p^{\alpha_i}\!-1)\ \ \Longrightarrow\ \
 \mathcal{V}_p(\underline{\alpha},D) = 0\,.\medskip \]
\item[c)] Keeping $\underline{\alpha}$ fixed, $\mathcal{V}_p(\underline{\alpha},D)$
    is monotonically decreasing in $D$.
\end{itemize}
\end{prop}

\noindent As we already determined
$\nu_p(\underline{\alpha},\underline{n})=\sum_{i=1}^{\smash{N}}\nu_p(\alpha_i,n_i)$ in
the lemma above, the precise calculation of $\mathcal{V}_p(\underline{\alpha},D)$ is
mere discrete optimization. We will do that in the next section, in Theorem
\ref{Thm.BC}. With the definitions of the numbers $\alpha:=\alpha_1+\dotsb+\alpha_N$
and $D_1,\dotsc,D_\alpha$ used in Theorem \ref{MAINTHM}, the result can be stated
as follows:
\begin{equation}
\mathcal{V}_p(\underline{\alpha},D)
\,=\,\alpha-\max\bigl\{0\leq t\leq \alpha\mid D_1+\dotsb+D_{t}\leq\tfrac{D}{p-1}\bigr\}\,.
\end{equation}

\subsection{The integral \protect{$\int_S f$}}

\noindent Let $A$ and $B$ be commutative groups, let $f \in B^A$\!, and let $S
\subseteq A$ be a finite subset. Following \cite{Karasev-Petrov12}, we set
\[ \int_S f \,\coloneqq\, \sum_{x \in S} f(x) \,\in\, B
 \quad\text{and}\quad \int\!f \,\coloneqq\, \int_A f\,. \]
Here we are mostly interested in the case $A=\Z^N$\!, $B=\Z$, and
$S=[p^{\underline{\alpha}})$. The following results generalize work of Wilson
\cite[Lemma 4]{Wilson06} about the case $\alpha_1=\dotsb=\alpha_N=1$\,:

\begin{prop}
\label{Thm.int} Let $D\in\N\cup\{\infty\}$ and $N,\alpha_1,\dotsc,\alpha_N\in \Z^+$\!. If
$f \in \Z^{\Z^N}$ has functional degree $\fdeg(f)\leq D$, then
\[ \ord_p\biggl(\int_{[p^{\underline{\alpha}})}f\biggr)
\,\geq\,\mathcal{V}_p\bigl(\underline{\alpha},D\bigr).
\]
\end{prop}
\begin{proof}
For commutative groups $A$ and $B$ and a finite subset $S \subseteq A$, the map
$\int_S: B^A \!\ra B$ is a $\Z$-module homomorphism - and this also holds when
$A=\Z^N$\!, $B = \Z$\,, and $S=[p^{\underline{\alpha}})$. By Theorem \ref{4.2}, it
therefore suffices to prove the inequality for functions of the form
\[ \underline{x} \,\mapsto\,\binom{x_1}{n_1} \dotsm \binom{x_N}{n_N} \]
with $|\underline{n}|\leq D$. This, however, is easy:
\[ \ord_p\biggl(\int_{[p^{\underline{\alpha}})}\!
 \binom{x_1}{n_1} \dotsm \binom{x_N}{n_N}\biggr)
=\,\prod_{j=1}^N \,\sum_{x_j \in [p^{\alpha_j})}\!\binom{x_j}{n_j}
=\,\nu_p\bigl(\underline{\alpha},\underline{n}\bigr)
\,\geq\,\mathcal{V}_p\bigl(\underline{\alpha},|\underline{n}|\bigr)
\,\geq\,\mathcal{V}_p\bigl(\underline{\alpha},D\bigr)\,.
\]
\end{proof}

\subsection{The Proof of Theorem \ref{MAINTHM}}\label{sec.proof}
Below is the proof of Theorem \ref{MAINTHM}, modulo two main discrete optimization
tasks. On one side, our proof shows that the broad outline of the argument is the same
as that of Theorem \ref{CS2MAIN2}, using the key ideas from Wilson's proof of Ax-Katz
over $\F_p$. On the other side, it motivates and sets up the new work of the present
paper, the two optimization tasks that are needed to complete the argument.

Our proof uses previously made definitions and the entire setup of Theorem
\ref{MAINTHM} without reintroducing them. Some new definitions are made on the way,
as well. As in \cite[\S6]{Aichinger-Moosbauer21}, we also use the tensor product of
functions: if $A_1,\ldots,A_n$ are commutative groups and $R$ is a rng, then the
\textbf{tensor product} $\bigotimes_{i=1}^n h_i$ of maps $h_i: A_i \ra R$ is the map
\[\bigotimes_{i=1}^n h_i:\,\bigoplus_{i=1}^n A_i \ra R\,,
 \,\  (x_1,x_2,\ldots,x_n) \mapsto h_1(x_1)h_2(x_2) \dotsm h_n(x_n)\,.\smallskip\]

\begin{proof}[Proof of Theorem \ref{MAINTHM}]
Let $\beta \in \Z^+$ be fixed given. For each $1 \leq j \leq r$ define the map $\chi_j: \Z
\ra \Z/p^{\beta} \Z$ by
\[ \chi_j(x) \,:=\, \begin{cases}
 1 & \text{if }x \equiv 0 \pmod{p^{\beta_j}}, \\ 0 & \text{otherwise,} \end{cases} \]
and let $\widetilde{\chi}_j: \Z \ra \Z$ be a proper lift of $\chi_j$ from $\Z/p^{\beta} \Z$ to
$\Z$\,. Using that $\Z/p^{\beta} \Z$ and $\Z$ are rngs, not just additive groups, set
$$\chi:=\,\bigotimes_{j=1}^r\chi_j\quad\text{and}\quad
 \widetilde{\chi}:=\,\bigotimes_{j=1}^r\widetilde{\chi}_j\,.$$
If $q$ denotes the quotient map from $\Z^N$ to
$A:=\bigoplus_{i=1}^N\Z/p^{\alpha_i}\Z$\,, and $\tilde{F}_j:\Z^N\!\ra\Z$ is a proper lift of
the pullback $F_j:\Z^N\!\ra\Z/p^{\beta_j}\Z$ of the function $f_j:A\ra\Z/p^{\beta_j}\Z$\,,
then
\[ \chi(\tilde{F}_1(\underline{x}),\dotsc,\tilde{F}_r(\underline{x}))
 \,=\, \begin{cases} 1 & \text{if}\ q(\underline{x}) \in Z(f_1,\dotsc,f_r), \\ 0
 & \text{otherwise,} \end{cases} \]
for each
$\underline{x}\in[p^{\underline{\alpha}}):=\prod_{i=1}^N\{0,1,\dotsc,p^{\alpha_i}\!-1\}$.\smallskip\\
Moreover, when restricted to $[p^{\underline{\alpha}})$, the quotient map $q$ induces a
bijection from $[p^{\underline{\alpha}})$ to $A$. Hence, with the function
$\tilde{\chi}(\tilde{F}_1,\dotsc,\tilde{F}_r):
 \underline{x}\mapsto\tilde{\chi}(\tilde{F}_1(\underline{x}),\dotsc,\tilde{F}_r(\underline{x}))$,
we get
\[\# Z(f_1,\dotsc,f_r)\,=\,kp^{\beta}+\int_{[p^{\underline{\alpha}})}
 \widetilde{\chi}(\tilde{F}_1,\dotsc,\tilde{F}_r)\quad\text{for some $k\in\Z$}\,.\]
We may certainly assume that $Z(f_1,\dotsc,f_r)$ is nonempty, so that $\ord_p \bigl(\#
Z(f_1,\dotsc,f_r)\bigr)$ is finite. Hence, after increasing our $\beta \in \Z^+$ if
necessary, we may assume
\begin{equation}\label{eq.beta}
\beta\,>\,\ord_p \bigl(\# Z(f_1,\dotsc,f_r)\bigr)\,\in\,\N\,,
\end{equation}
and with that
\begin{equation*}
\ord_p \bigl(\# Z(f_1,\dotsc,f_r)\bigr)\,=\,
\ord_p \biggl(\int_{[p^{\underline{\alpha}})} \widetilde{\chi}(\tilde{F}_1,\dotsc,\tilde{F}_r)\biggr)\,.
\end{equation*}
Now, for each $1 \leq j \leq r$, Corollary \ref{COR4.57} provides an integer valued
function $c_j$ on the set
$$[\hat n_j(\beta)]\,:=\,\{0,1,\dotsc,\hat n_j(\beta)\}
 \quad\text{where}\quad\hat n_j(\beta) \,:=\,(p^{\beta_j}\!-1)+(\beta-1)p^{\beta_j-1}(p-1)\,,$$
such that, for each $x\in\Z$\,,
\[\tilde{\chi}_j(x) \,=\! \sum_{n \in [\hat n_j(\beta)]} \binom{x}{n} c_j(n)\,. \]
With
$$\hat{\underline{n}}(\beta):=(\hat n_1(\beta),\ldots,\hat n_r(\beta))
 \quad\text{and}\quad[\hat{\underline{n}}(\beta)]:=\prod_{j=1}^r[\hat n_j(\beta)]
 \,\subseteq\,\N^r\!,$$
for each $\underline{x}\in\Z^N$\!,
\begin{eqnarray*}
\tilde{\chi}(\tilde{F}_1(\underline{x}),\dotsc,\tilde{F}_r(\underline{x})) &={}&
\tilde{\chi}_1(\tilde{F_1}(\underline{x})) \dotsm \tilde{\chi}_r(\tilde{F_r}(\underline{x})) \\
&={}&\!\sum_{\nn \in[\underline{\hat
n}(\beta)]} \binom{ \tilde{F}_1(\underline{x})}{n_1} \dotsm
\binom{\tilde{F}_r(\underline{x})}{n_r}c_1(n_1) \dotsm c_r(n_r)\,.
\end{eqnarray*}
Hence, with the functions $\binom{ \tilde{F}_j}{n_j}:
 \underline{x}\mapsto\binom{ \tilde{F}_j(\underline{x})}{n_j}$,
\begin{eqnarray*}
\int_{[p^{\underline{\alpha}})} \tilde{\chi}(\tilde{F}_1,\dotsc,\tilde{F}_r)
&={}&\!\sum_{\nn \in[\underline{\hat n}(\beta)]} c_1(n_1) \dotsm c_r(n_r)
\int_{[p^{\underline{\alpha}})}\binom{ \tilde{F}_1}{n_1} \dotsm \binom{\tilde{F}_r}{n_r}\,.
\end{eqnarray*}
So if we put
$$\mathbf{m} \,:=\,\min_{\underline{n}\in[\underline{\hat
n}(\beta)]}\biggl(\ord_p \bigl(c_1(n_1)\bigr) + \dotsb +
\ord_p \bigl(c_r(n_r)\bigr) \,+\, \ord_p \biggl( \int_{[p^{\underline{\alpha}})}
\binom{ \tilde{F}_1}{n_1} \dotsm \binom{\tilde{F}_r}{n_r}\biggr)\biggr),
$$
it follows that
$$\ord_p \bigl(\# Z(f_1,\dotsc,f_r)\bigr)\,=\, \ord_p \biggl(\int_{[p^{\underline{\alpha}})}
 \widetilde{\chi}(\tilde{F}_1,\dotsc,\tilde{F}_r)\biggr) \,\geq\, \mathbf{m}\,.$$
Thus the matter of it is to give a good lower bound on the quantity $\mathbf{m}$,
using that $\fdeg(\tilde{F}_j) = \fdeg(f_j) \leq d_j$ for all $1 \leq j \leq r$ (cf.\ \cite[Cor.\
2.13, \S 2.4 and \S 2.5]{Clark-Schauz23a}).  Part of this can be quickly done in the
same way as in \cite{Clark-Schauz23a}: Corollary \ref{COR4.57} also says that the
functions $c_j: [\hat n_j(\beta)] \ra \Z$ can be chosen such that, for each $h \in \Z^+$\!
and $n\in[\hat n_j(\beta)]$,
\begin{equation*}
\ p^{\beta_j}\!-1+p^{\beta_j-1}(p-1)(h-1)<n
\,\implies\, p^h \bigm{|} c_j(n)\,.
\end{equation*}
Taking
$$h_j\,=\,h_j(n_j)\,:=\,\biggl\lceil\frac{n_j-(p^{\beta_j}\!-1)}{p^{\beta_j-1}(p-1)}\biggr\rceil
\,<\,\frac{n_j-(p^{\beta_j}\!-1)}{p^{\beta_j-1}(p-1)}+1,$$ we have
\[ p^{\beta_j}\!-1 + (h_j-1)p^{\beta_j-1}(p-1) \,<\, n_j\,, \]
and thus Corollary \ref{COR4.57} yields
\begin{equation*}
\ord_p \bigl(c_j(n_j)\bigr)\,\geq\, \overbar{h_j}
\,=\,\overbar{\biggl\lceil\frac{n_j-(p^{\beta_j}\!-1)}{p^{\beta_j-1}(p-1)}\biggr\rceil},
\end{equation*}
where, for real numbers $h$,
\begin{equation*}
\overbar{h}\,:=\,\max(h,0)\,.
\end{equation*}
\\
Moreover, using
\cite[Thm.\ 4.3 and Lem.\ 6.1]{Aichinger-Moosbauer21}, we have
\begin{align*}
\fdeg\biggl( \binom{\tilde{F}_1}{n_1} \dotsm \binom{\tilde{F}_r}{n_r}\biggr)
&\,\leq\,\,\sum_{j=1}^rd_jn_j\,,
\end{align*}
and Proposition \ref{Thm.int} shows that
\begin{equation*}
\ord_p\biggl( \int_{[p^{\underline{\alpha}})} \binom{\tilde{F_1}}{n_1} \dotsm \binom{\tilde{F_r}}{n_r} \biggr)
\,\geq\,\,\mathcal{V}_p\bigl(\underline{\alpha},{\textstyle\sum_{j=1}^r}d_jn_j\bigr)\,.
\end{equation*}
We deduce that
\begin{align*}
\ord_p(\#Z(f_1,\ldots,f_r)) \,\geq\, \mathbf{m}
\,\geq\,\min_{\underline{n}\in[\underline{\hat n}(\beta)]}\mathcal N(\underline{n})
\end{align*}
where
$$\mathcal N(\underline{n}) \,:=\,
\sum_{j=1}^r\overbar{\biggl\lceil\frac{n_j-(p^{\beta_j}\!-1)}{p^{\beta_j-1}(p-1)}\biggr\rceil}
 +\, \mathcal{V}_p\Bigl(\underline{\alpha},\,{\textstyle\sum_{j=1}^r}d_jn_j\Bigr).$$
The precise calculation of the minimum of $\mathcal N(\underline{n})$ when
$\underline{n}$ runes through $[\underline{\hat n}(\beta)]$ is mere discrete
optimization. We will do that in Section 5, in Lemma \ref{Lem.S}. After increasing
$\beta$ if necessary\footnote{We need
$\beta>s_0:=\overbar{\bigl\lceil(\breve{\mathcal{A}}-\mathcal B)/
 (d_1p^{{\beta_1-1}})\bigr\rceil}$ in Lemma \ref{Lem.S}.
Within the full proof of Theorem \ref{MAINTHM}, however, we assume
$\beta>\ord_p(\# Z)$ already in \eqref{eq.beta}, and $\ord_p(\# Z)\geq s_0$ by the
findings of this paper.}\!, it yields
$$
\min_{\underline{n}\in[\underline{\hat n}(\beta)]}\mathcal N(\underline{n})
 \,=\,\begin{cases}
 \displaystyle{\biggl\lceil\frac{\breve{\mathcal{A}}-\mathcal B}{d_1p^{{\beta_1-1}}}\biggr\rceil
 +\alpha-\breve\alpha}
 & \text{if $\breve{\mathcal{A}}>\mathcal B$,}\\
 \,\alpha-\max\bigl\{1\leq t\leq\alpha\mid D_1+\dotsb+D_{t}\leq \mathcal B\bigr\}^{\strut}\!
 & \text{if $\breve{\mathcal{A}}\leq\mathcal B$.}
\end{cases}
$$
Notice that the answer obtained is independent of $\beta$.
\end{proof}

\subsection{The Proof of Corollary \ref{MAINCOR}}\label{sec.proof2}
As claimed before, the result takes a somewhat simpler form when $\alpha_1 =
\dotsb = \alpha_N$.

\begin{proof}[Proof of Corollary \ref{MAINCOR}]
If $\alpha_1 = \dotsb = \alpha_N$ then
$$\alpha\,=\,N\alpha_1\,,\quad\breve\alpha\,=\,N\breve\alpha_1\,,\quad
 \alpha_1' \,=\,\dotsb\,=\,\alpha_{\alpha_1}'\,=\,N\quad\text{and}\quad
 \breve{\mathcal{A}}\,=\,N\frac{p^{\breve\alpha_1}-1}{p-1}\,,$$
which yields the claimed simplifications in the case
$\breve{\mathcal{A}}>\mathcal B$. With the parameters
$$Q\,:=\,\bigl\lfloor\log_p\bigl((p-1)\mathcal B/N+1\bigr)\bigr\rfloor\quad\text{and}\quad
 R\,:=\,\biggl\lfloor\frac{\mathcal B-N\frac{p^Q-1}{p-1}}{p^Q}\biggr\rfloor$$
we also have
$$Q\,\leq\,\log_p\bigl((p-1)\mathcal B/N+1\bigr)\,<\,Q+1$$
i.e.,$$ N\frac{p^Q-1}{p-1}\,\leq\,\mathcal B
 \,<\,N\frac{p^{Q+1}-1}{p-1}=N\frac{p^Q-1}{p-1}+Np^Q
$$
and
$$R\,\leq\,\frac{\mathcal B-N\frac{p^Q-1}{p-1}}{p^Q}\,<\,R+1$$
i.e.,
\begin{equation}\label{eq.QR}
N\frac{p^Q-1}{p-1}+Rp^Q\,\leq\,\mathcal B
 \,<\,N\frac{p^Q-1}{p-1}+(R+1)p^Q\,.
\end{equation}
In particular, $N\frac{p^Q-1}{p-1}\leq\mathcal B$ and
$N\frac{p^Q-1}{p-1}+Rp^Q\leq\mathcal B<N\frac{p^Q-1}{p-1}+Np^Q$, so that, on one
side,
$$0\,\leq\,R\,<\,N\,.$$
If we further assume $D_1+\dotsb+D_{\alpha}>\mathcal B$, then
$N\frac{p^{\alpha_1}-1}{p-1}=D_1+\dotsb+D_{\alpha}>\mathcal B\geq
N\frac{p^Q-1}{p-1}$, so that, on the other side,
$$0\,\leq\,Q\,<\,\alpha_1\,.$$
From this follows $QN+R<N\alpha_1=\alpha$, and we see that
$D_{QN+1}=\dotsb=D_{QN+R+1}=p^Q$, i.e.,
$D_1+\dotsb+D_{QN+R}=N\frac{p^Q-1}{p-1}+Rp^Q$ and
$D_1+\dotsb+D_{QN+R+1}=N\frac{p^Q-1}{p-1}+(R+1)p^Q$. Hence, Inequality
\eqref{eq.QR} can be restated as
$$D_1+\dotsb+D_{QN+R}\,\leq\,\mathcal B\,<\,D_1+\dotsb+D_{QN+R+1}\,,$$
which means that
$$\max\bigl\{1\leq t\leq\alpha\mid D_1+\dotsb+D_{t}\leq \mathcal B\bigr\}\,=\,NQ+R\,.$$
Thus, we can replace the second lower bound $\alpha-\max\bigl\{1\leq t\leq\alpha\mid
D_1+\dotsb+D_{t}\leq \mathcal B\bigr\}$ with $N(\alpha_1-Q)-R$\,. This replacement is
also correct in the case $D_1+\dotsb+D_{\alpha}\leq\mathcal B$\,, because then both
terms are non-positive. This is clear for $\alpha-\max\bigl\{1\leq t\leq\alpha\mid
D_1+\dotsb+D_{t}\leq \mathcal B\bigr\}$, but we also see that
$N\frac{p^{\alpha_1}-1}{p-1}=D_1+\dotsb+D_{\alpha}\leq\mathcal
B<N\frac{p^{Q+1}-1}{p-1}$ implies $\alpha_1<Q+1$, which entails
$N(\alpha_1-Q)-R\leq0$\,. We obtain
$$
\ord_p \bigl(\# Z_{A}(f_1,\dotsc,f_r))
 \,\geq\,\begin{cases}
 \displaystyle{\biggl\lceil\frac{N\frac{p^{\breve\alpha_1}-1}{p-1}-\mathcal B}{d_1p^{{\beta_1-1}}}\biggr\rceil}
 \,+\,N(\alpha_1-\breve\alpha_1)
 & \text{if $N\frac{p^{\breve\alpha_1}-1}{p-1}>\mathcal B$,}\\
 \,N(\alpha_1-Q)-R^{\strut\phantom{|}}
 & \text{if $N\frac{p^{\breve\alpha_1}-1}{p-1}\leq\mathcal B$.}
\end{cases}
$$
\end{proof}

\section{Minimization of \protect{$\nu_p(\protect\underline{\alpha},\bullet)$}}

\noindent In this section, we determine the minimum value
$\mathcal{V}_p(\underline{\alpha},D)$ of the function
$\nu_p(\underline{\alpha},\bullet)$ over the restricted domain
$$\mathcal{D}(N,D)\,:=\,\bigl\{\underline{n}\in\N^N\!\bigm{|}|\underline{n}|\leq D\bigr\}\,,$$
where the numbers $N,\alpha_1,\dotsc,\alpha_N\in \Z^+$ with
$\alpha_1\geq\dotsb\geq\alpha_N$ and $D\in\N$ are fixed given. (The case $D=\infty$
is trivial, as $\mathcal{V}_p(\underline{\alpha},\infty)=0$ by Proposition \ref{Pr.Vv}\,b.)
In our investigation, the original definition of $\nu_p(\underline{\alpha},\bullet)$ does not
actually matter. We may view the formula in Lemma \ref{Lem.vp} as the definition. More
precisely, for $\underline{n}\in\N^N$\!,
$$\nu_p(\underline{\alpha},\underline{n})\,:=\,\sum_{i=1}^{N}\nu_p(\alpha_i,n_i)
\quad\text{with}\quad\nu_p(\alpha_i,n_i)
\,:=\,\begin{cases}
\alpha_i-\ord_p(n_i+1) & \text{if $n_i\leq p^{\alpha_i}\!-1$,}\\
\infty & \text{otherwise.}
\end{cases}
$$
Our final result will be stated in terms of the parameters $\alpha:=\alpha_1+ \dotsb +
\alpha_N$, and $D_1,\dotsc,D_\alpha$ of Theorem \ref{MAINTHM}, i.e.,
$$\bigl(D_1,D_2,\dotsc,D_{\alpha}\bigr)
\,:=\,\bigl(\,\underbrace{1,1,\dotsc,1}_{\alpha_1'\ \text{times}},
\,\underbrace{p,p,\dotsc,p}_{\alpha_2'\ \text{times}},
\,\dotsc,\,\underbrace{p^{\alpha_1-1},p^{\alpha_1-1},\dotsc,p^{\alpha_1-1}}%
 _{\alpha_{\alpha_1}'\ \text{times}}\,\bigr).$$

\subsection{The Minimum Value \protect{$\mathcal{V}_p(\protect\underline{\alpha},D)$} of
\protect{$\nu_p(\protect\underline{\alpha},\bullet)$} over \protect{$\mathcal{D}(N,D)$}}

\begin{thm}
\label{Thm.BC} In the setting above, with $D\in\N$\,, the function
$$\nu_p(\underline{\alpha},\bullet)\bigr{|}_{\mathcal{D}(N,D)}
 \!:\,\mathcal{D}(N,D) \,\longrightarrow\,\N\cup\{\infty\}\,,\ \ \underline{n}
 \,\longmapsto\,\nu_p(\underline{\alpha},\underline{n})$$
has minimum value
$$\mathcal{V}_p(\underline{\alpha},D)
 \,=\,\alpha-\max\bigl\{0\leq t\leq\alpha\mid D_1+\dotsb+D_{t}\leq\tfrac{D}{p-1}\bigr\}\,.$$
\end{thm}

\begin{proof}We may restrict the domain of $\nu_p(\underline{\alpha},\bullet)$ from
$\mathcal{D}(N,D)$ to $\mathcal{D}(N,D)\cap[p^{\underline{\alpha}})$ with
$[p^{\underline{\alpha}}):
 =\prod_{i=1}^N\{0,1,\dotsc,p^{\alpha_i}\!-1\}$,
because $\nu_p(\underline{\alpha},\bullet)$ is finite inside but positive infinite outside of
$[p^{\underline{\alpha}})$. Inside $[p^{\underline{\alpha}})$, however,
$$\nu_p(\underline{\alpha},\underline{n})=\alpha-\sum_{i=1}^{N}\ord_p(n_i+1)\,.$$
So, we need to find the maximum of the function
$$\lambda_p:\mathcal{D}(N,D)\cap[p^{\underline{\alpha}})\longrightarrow\N\,,\quad
 \underline{n}\longmapsto\lambda_p(\underline{n}):=\sum_{i=1}^{N}\ord_p(n_i+1)\,.$$

If the point $\underline{n}=(n_i)_{i=1}^N$ of the domain
$\mathcal{D}(N,D)\cap[p^{\underline{\alpha}})$ happens to be a maximum point, then
the point $\underline{\tilde n}=(\tilde n_i)_{i=1}^N$ with $\tilde
n_i+1:=p^{\ord_p(n_i+1)}$ is also a maximum point in that domain, because
$\ord_p(\tilde n_i+1)=\ord_p(n_i+1)$ and $0\leq\tilde n_i\leq n_i$ for all $1\leq i\leq N$.
Hence, we may restrict our attention to points $\underline{n}$ with the property that
each $n_i+1$ is a power of $p$, say $n_i+1=p^{\mu_i}$. With the substitutions
$n_i:=p^{\mu_i}-1$ in mind, we then just have to find the maximum of the function
$$\lambda:\bigl\{\underline{\mu}\in[\underline{\alpha}]\bigm|
 \omega(\underline{\mu})\leq\tfrac{D}{p-1}\bigr\}\longrightarrow\N\ ,\quad
 \underline{\mu}\longmapsto\lambda(\underline{\mu}):=\lambda_p\bigl((p^{\mu_i}-1)_{i=1}^N\bigr)
 =\sum_{i=1}^{N}\mu_i\,,$$
where $[\underline{\alpha}]:=\prod_{i=1}^N\{0,1,\dotsc,\alpha_i\}$ and
$$\omega(\underline{\mu})\,:=\,\sum_{i=1}^{N}\frac{p^{\mu_i}-1}{p-1}
 \,=\,\sum_{i=1}^N\sum_{j=0}^{\mu_i-1}p^j\,=\,\sum_{i=1}^N\sum_{j=1}^{\mu_i}p^{j-1}.$$

Now, if we draw Ferrers-type diagrams for the potential arguments $\underline{\mu}$ of
$\lambda$ as sub-diagrams of Ferrers' diagram of $\underline{\alpha}$ (representing
$\mu_i$ by $\mu_i$ consecutive dots in row $i$) then $\lambda(\underline{\mu})$ is the
number of dots in the sub-diagrams of $\underline{\mu}$, while
$\omega(\underline{\mu})$ gives a weighted count of those dots -- a dot in the
$j^\text{th}$ column is counted with weight $p^{j-1}$\!, as shown in Figure \ref{fig.w}.
Hence, to find the maximum of $\lambda$, we need to maximize the number of dots in
the sub-diagram corresponding to $\underline{\mu}$\,, while keeping their total weight
(total cost) $\omega(\underline{\mu})$ below $\frac{D}{p-1}$\,. So, when selecting the
dots in $\underline{\mu}$\,, we just have to select the cheapest dots first. In our situation
of column-wise increasing weights, we have to select the dots column by column, from
left to right, starting with the left-most column with dots of lowest weight. Insight a
column the order of selection does not matter, as long as the column is completely
finished before we move to the next column. We may just go top-down inside columns,
as in Figure \ref{fig.w}. Following that order, we collect in step $t$ a dot of weight
$D_t$\,, because that is how we defined $D_t$\,. Hence, after $t$ steps we obtain a
$\underline{\mu}$ with
$$\omega(\underline{\mu})=D_1+\dotsb+D_t
 \quad\text{and}\quad\lambda(\underline{\mu})=t\,.$$
Our selection process has to stop when the limit $\frac{D}{p-1}$ for
$\omega(\underline{\mu})$ is reached, that is when
$$t\,=\,t(D)\,:=\,\max\bigl\{0\leq t\leq\alpha\mid D_1+\dotsb+D_{t}\leq\tfrac{D}{p-1}\bigr\}\,.$$
At that point, $\underline{\mu}=\underline{\mu}(D)$ is a maximum point of $\lambda$,
and the associated $\underline{n}(D):=(p^{\mu_i}-1)_{i=1}^N$ is a minimum point of
$\nu_p(\underline{\alpha},\bullet)$ in $\mathcal{D}(N,D)$. The minimum value is
$$\mathcal{V}_p(\underline{\alpha},D)\,=\,\nu_p(\underline{\alpha},\underline{n}(D))
 \,=\,\alpha-\lambda(\underline{\mu})
 \,=\,\alpha-\max\bigl\{0\leq t\leq\alpha\mid D_1+\dotsb+D_{t}\leq\tfrac{D}{p-1}\bigr\}\,.$$
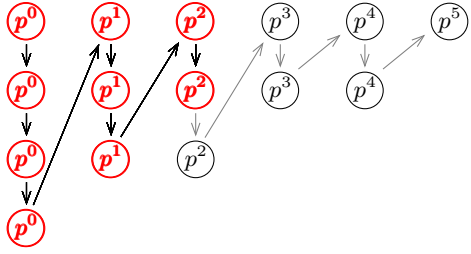
\begin{figure}[h!]
\begin{center}
\[\small\begin{tikzcd}[sep=10pt, row sep=5pt]
        \pmb{\color{red}\cp{0}\ar[d, shorten <=-2pt, shorten >=-2pt, -{Straight Barb[scale length=2, scale width=.7]}]} & \pmb{\color{red}\cp{1}\ar[d, shorten <=-2pt, shorten >=-2pt, -{Straight Barb[scale length=2, scale width=.7]}]} & \pmb{\color{red}\cp{2}\ar[d, shorten <=-2pt, shorten >=-2pt, -{Straight Barb[scale length=2, scale width=.7]}]} & \cp{3}\ar[d, color=gray, shorten <=-2pt, shorten >=-2pt, -{Straight Barb[scale length=2, scale width=.7]}] & \cp{4}\ar[d, color=gray, shorten <=-2pt, shorten >=-2pt, -{Straight Barb[scale length=2, scale width=.7]}] & \cp{5}\\
        \pmb{\color{red}\cp{0}\ar[d, shorten <=-2pt, shorten >=-2pt, -{Straight Barb[scale length=2, scale width=.7]}]} & \pmb{\color{red}\cp{1}\ar[d, shorten <=-2pt, shorten >=-2pt, -{Straight Barb[scale length=2, scale width=.7]}]} & \pmb{\color{red}\cp{2}}\ar[d, color=gray, shorten <=-2pt, shorten >=-2pt, -{Straight Barb[scale length=2, scale width=.7]}] & \cp{3}\ar[ur, color=gray, shorten <=-4pt, shorten >=-2.6pt, -{Straight Barb[scale length=2, scale width=.7]}] & \cp{4}\ar[ur, color=gray, shorten <=-4pt, shorten >=-2.6pt, -{Straight Barb[scale length=2, scale width=.7]}] &\\
        \pmb{\color{red}\cp{0}\ar[d, shorten <=-2pt, shorten >=-2pt, -{Straight Barb[scale length=2, scale width=.7]}]} & \pmb{\color{red}\cp{1}\ar[uur, shorten <=-2.8pt, shorten >=-3.3pt, -{Straight Barb[scale length=2, scale width=.7]}]} & \cp{2}\ar[uur, color=gray, shorten <=-2.8pt, shorten >=-3.3pt, -{Straight Barb[scale length=2, scale width=.7]}] &&\\
        \pmb{\color{red}\cp{0}\ar[uuur, shorten <=-2.2pt, shorten >=-2.2pt, -{Straight Barb[scale length=2, scale width=.7]}]} &&&&
\end{tikzcd}\vspace{-1em}\]
\end{center}
\caption{\label{fig.w}The minimum weight of a set of $9$ dots inside $\underline{\alpha}=(6,5,3,1)$ is $D_1+D_2+\dotsb+D_9=4+3p+2p^2$.}
\end{figure}
\end{proof}

\noindent Using our new formula, we can now show that
$D<\sum_{i=1}^N\frac{p^{\alpha_i}\!-1}{p-1}$ is not just necessary for
$\mathcal{V}_p(\underline{\alpha},D)>0$, as we already have seen in Proposition
\ref{Pr.Vv}\,b, it is also sufficient:

\begin{cor}\label{Cor.eq}
Maintain the setup of Theorem \ref{Thm.BC}, we have
$$
\mathcal{V}_p(\underline{\alpha},D)>0\quad\Longleftrightarrow\quad
D<\sum_{i=1}^N\bigl(p^{\alpha_i}\!-1\bigr)\,.$$
\end{cor}

\begin{proof}
With the last result in Example \ref{Ex.con1}, we see that
$$D_1+\dotsb+D_{\alpha}
\,=\,\alpha_1'p^0+\dotsb+\alpha_{N}'p^{N-1}
 \,=\,\sum_{i=1}^N\frac{p^{\alpha_i}\!-1}{p-1}\,.$$
So, by Theorem \ref{Thm.BC},
\begin{align*}
\mathcal{V}_p(\underline{\alpha},D)>0
 \ \,&\Longleftrightarrow\ \,
 \max\bigl\{0\leq t\leq\alpha\mid D_1+\dotsb+D_{t}\leq\tfrac{D}{p-1}\bigr\}<\alpha\\
 \ \,&\Longleftrightarrow\ \,
 D_1+\dotsb+D_{\alpha}\nleq\tfrac{D}{p-1}\\
 \ \,&\Longleftrightarrow\ \,
 D<\sum_{i=1}^N\bigl(p^{\alpha_i}\!-1\bigr)\,.\qedhere
\end{align*}
\end{proof}

\subsection{Alternative Expressions for \protect{$\mathcal{V}_p(\protect\underline{\alpha},D)$}
and Special Cases.} 
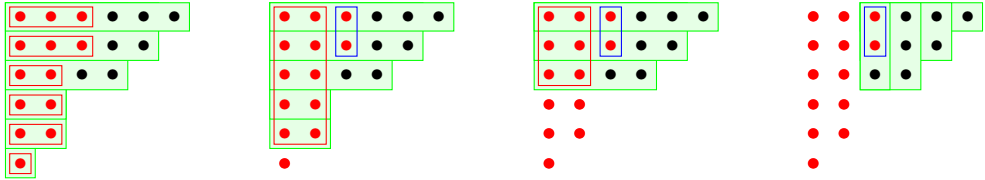
\begin{figure}[h!]
\begin{center}
\begin{tikzpicture}[baseline=(current  bounding  box.center)]
\matrix [matrix of math nodes] (m)
{
        {\color{red}\bullet} & {\color{red}\bullet} & {\color{red}\bullet} & \bullet & \bullet & \bullet\\
        {\color{red}\bullet} & {\color{red}\bullet} & {\color{red}\bullet} & \bullet & \bullet & \\
        {\color{red}\bullet} & {\color{red}\bullet} & \bullet & \bullet & & \\
        {\color{red}\bullet} & {\color{red}\bullet} & & & & \\
        {\color{red}\bullet} & {\color{red}\bullet} & & & & \\
        {\color{red}\bullet} & & & & &\\
};
\node[draw,fit=(m-1-1) (m-6-1),green,fill=green!10,inner sep=-0.4pt]{};
\node[draw,fit=(m-1-1) (m-5-2),green,fill=green!10,inner sep=-0.4pt]{};
\node[draw,fit=(m-1-1) (m-4-2),green,fill=green!10,inner sep=-0.4pt]{};
\node[draw,fit=(m-1-1) (m-3-4),green,fill=green!10,inner sep=-0.4pt]{};
\node[draw,fit=(m-1-1) (m-2-5),green,fill=green!10,inner sep=-0.4pt]{};
\node[draw,fit=(m-1-1) (m-1-6),green,fill=green!10,inner sep=-0.4pt]{};
\node[draw,fit=(m-1-1) (m-1-3),red,inner sep=-2pt]{};
\node[draw,fit=(m-2-1) (m-2-3),red,inner sep=-2pt]{};
\node[draw,fit=(m-3-1) (m-3-2),red,inner sep=-2pt]{};
\node[draw,fit=(m-4-1) (m-4-2),red,inner sep=-2pt]{};
\node[draw,fit=(m-5-1) (m-5-2),red,inner sep=-2pt]{};
\node[draw,fit=(m-6-1) (m-6-1),red,inner sep=-2pt]{};
\matrix [matrix of math nodes] (m)
{
        {\color{red}\bullet} & {\color{red}\bullet} & {\color{red}\bullet} & \bullet & \bullet & \bullet\\
        {\color{red}\bullet} & {\color{red}\bullet} & {\color{red}\bullet} & \bullet & \bullet & \\
        {\color{red}\bullet} & {\color{red}\bullet} & \bullet & \bullet & & \\
        {\color{red}\bullet} & {\color{red}\bullet} & & & & \\
        {\color{red}\bullet} & {\color{red}\bullet} & & & & \\
        {\color{red}\bullet} & & & & &\\
};
\end{tikzpicture}
\qquad
\begin{tikzpicture}[baseline=(current  bounding  box.center)]
\matrix [matrix of math nodes] (m)
{
        {\color{red}\bullet} & {\color{red}\bullet} & {\color{red}\bullet} & \bullet & \bullet & \bullet\\
        {\color{red}\bullet} & {\color{red}\bullet} & {\color{red}\bullet} & \bullet & \bullet & \\
        {\color{red}\bullet} & {\color{red}\bullet} & \bullet & \bullet & & \\
        {\color{red}\bullet} & {\color{red}\bullet} & & & & \\
        {\color{red}\bullet} & {\color{red}\bullet} & & & & \\
        {\color{red}\bullet} & & & & &\\
};
\node[draw,fit=(m-1-1) (m-5-2),green,fill=green!10,inner sep=-0.4pt]{};
\node[draw,fit=(m-1-1) (m-4-2),green,fill=green!10,inner sep=-0.4pt]{};
\node[draw,fit=(m-1-1) (m-3-4),green,fill=green!10,inner sep=-0.4pt]{};
\node[draw,fit=(m-1-1) (m-2-5),green,fill=green!10,inner sep=-0.4pt]{};
\node[draw,fit=(m-1-1) (m-1-6),green,fill=green!10,inner sep=-0.4pt]{};
\node[draw,fit=(m-1-1) (m-5-2),red,inner sep=-2pt]{};
\node[draw,fit=(m-1-3) (m-2-3),blue,inner sep=-2pt]{};
\matrix [matrix of math nodes] (m)
{
        {\color{red}\bullet} & {\color{red}\bullet} & {\color{red}\bullet} & \bullet & \bullet & \bullet\\
        {\color{red}\bullet} & {\color{red}\bullet} & {\color{red}\bullet} & \bullet & \bullet & \\
        {\color{red}\bullet} & {\color{red}\bullet} & \bullet & \bullet & & \\
        {\color{red}\bullet} & {\color{red}\bullet} & & & & \\
        {\color{red}\bullet} & {\color{red}\bullet} & & & & \\
        {\color{red}\bullet} & & & & &\\
};
\end{tikzpicture}
\qquad
\begin{tikzpicture}[baseline=(current  bounding  box.center)]
\matrix [matrix of math nodes] (m)
{
        {\color{red}\bullet} & {\color{red}\bullet} & {\color{red}\bullet} & \bullet & \bullet & \bullet\\
        {\color{red}\bullet} & {\color{red}\bullet} & {\color{red}\bullet} & \bullet & \bullet & \\
        {\color{red}\bullet} & {\color{red}\bullet} & \bullet & \bullet & & \\
        {\color{red}\bullet} & {\color{red}\bullet} & & & & \\
        {\color{red}\bullet} & {\color{red}\bullet} & & & & \\
        {\color{red}\bullet} & & & & &\\
};
\node[draw,fit=(m-1-1) (m-3-4),green,fill=green!10,inner sep=-0.4pt]{};
\node[draw,fit=(m-1-1) (m-2-5),green,fill=green!10,inner sep=-0.4pt]{};
\node[draw,fit=(m-1-1) (m-1-6),green,fill=green!10,inner sep=-0.4pt]{};
\node[draw,fit=(m-1-1) (m-3-2),red,inner sep=-2pt]{};
\node[draw,fit=(m-1-3) (m-2-3),blue,inner sep=-2pt]{};
\matrix [matrix of math nodes] (m)
{
        {\color{red}\bullet} & {\color{red}\bullet} & {\color{red}\bullet} & \bullet & \bullet & \bullet\\
        {\color{red}\bullet} & {\color{red}\bullet} & {\color{red}\bullet} & \bullet & \bullet & \\
        {\color{red}\bullet} & {\color{red}\bullet} & \bullet & \bullet & & \\
        {\color{red}\bullet} & {\color{red}\bullet} & & & & \\
        {\color{red}\bullet} & {\color{red}\bullet} & & & & \\
        {\color{red}\bullet} & & & & &\\
};
\end{tikzpicture}
\qquad
\begin{tikzpicture}[baseline=(current  bounding  box.center)]
\matrix [matrix of math nodes] (m)
{
        {\color{red}\bullet} & {\color{red}\bullet} & {\color{red}\bullet} & \bullet & \bullet & \bullet\\
        {\color{red}\bullet} & {\color{red}\bullet} & {\color{red}\bullet} & \bullet & \bullet & \\
        {\color{red}\bullet} & {\color{red}\bullet} & \bullet & \bullet & & \\
        {\color{red}\bullet} & {\color{red}\bullet} & & & & \\
        {\color{red}\bullet} & {\color{red}\bullet} & & & & \\
        {\color{red}\bullet} & & & & &\\
};
\node[draw,fit=(m-1-3) (m-1-6),green,fill=green!10,inner sep=-0.4pt]{};
\node[draw,fit=(m-1-3) (m-2-5),green,fill=green!10,inner sep=-0.4pt]{};
\node[draw,fit=(m-1-3) (m-3-4),green,fill=green!10,inner sep=-0.4pt]{};
\node[draw,fit=(m-1-3) (m-3-3),green,fill=green!10,inner sep=-0.4pt]{};
\node[draw,fit=(m-1-3) (m-2-3),blue,inner sep=-2pt]{};
\matrix [matrix of math nodes] (m)
{
        {\color{red}\bullet} & {\color{red}\bullet} & {\color{red}\bullet} & \bullet & \bullet & \bullet\\
        {\color{red}\bullet} & {\color{red}\bullet} & {\color{red}\bullet} & \bullet & \bullet & \\
        {\color{red}\bullet} & {\color{red}\bullet} & \bullet & \bullet & & \\
        {\color{red}\bullet} & {\color{red}\bullet} & & & & \\
        {\color{red}\bullet} & {\color{red}\bullet} & & & & \\
        {\color{red}\bullet} & & & & &\\
};
\end{tikzpicture}
\end{center}
\caption{\label{fig.w2}The four ways to calculate $\mathcal{V}_p(\underline{\alpha},D)$
(black dots) in Rem.\ \ref{rem.V2}. Here $\underline{\alpha}=(6,5,4,2,2,1)$ and
$D$ is such that $t(D)=13$ (red dots). So, $\color{red} \underline{\mu}(D)=(3,3,2,2,2,1)$,
$\color{red} Q(D)=2$, $\color{red} \alpha'_{Q(D)}=5$, $\color{red} \alpha'_{Q(D)+1}=3$,
$\color{blue} R(D)=2$.}
\end{figure}

In our main theorem (Theorem \ref{MAINTHM}), we presented only
one formula as the final result, as we did not want to make things any more complicated
than necessary. That result can, however, be stated in different forms, by replacing the
second lower bound \[\alpha-\max\bigl\{0\leq t\leq\alpha\mid
D_1+\dotsb+D_{t}\leq\tfrac{D}{p-1}\bigr\}
 =\mathcal{V}_p(\underline{\alpha},D) \]
with an alternative expression for $\mathcal{V}_p(\underline{\alpha},D)$. We present
several alternative formulas in the second remark below (and in the subsequent
corollary), after extracting an additional insights from our previous calculation of
$\mathcal{V}_p(\underline{\alpha},D)$ in the following first remark:

\begin{remark} \label{rem.V}\mbox{}
In the proof of Theorem \ref{Thm.BC} we also constructed a minimum point
$\underline{n}(D)$ of the function $\mathcal{D}(N,D)\ra\N\cup\{\infty\}$\,,
$\underline{n}\mapsto\nu_p(\underline{\alpha},\underline{n})$\,. This point may be
written as $$\underline{n}(D)\,:=\, \bigl(p^{\mu_i(D)}\!-1\bigr)_{i=1}^N\,,$$
where\vspace{-.2em}
\begin{equation*}
\mu_i(D)\,:=\,\begin{cases}
Q(D)+1 & \text{if \,$1\leq i\leq R(D)$}\\
Q(D) & \text{if \,$R(D)<i\leq\alpha_{Q(D)+1}'$}\\
\,\alpha_i & \text{if \,$\alpha_{Q(D)+1}'\!< i\leq N$}
\end{cases}
\end{equation*}
with\vspace{-.2em}
\begin{align*}
Q(D)\,:={}&\,\max\Bigl\{\textstyle0\leq Q\leq\alpha_1\bigm{|}
D_1+\dotsb+D_{\alpha_1'+\dotsb+\alpha_{Q}'}\leq \frac{D}{p-1}\Bigr\}\\
={}&\,\max\Bigl\{\textstyle0\leq Q\leq\alpha_1\bigm{|}
\sum_{j=1}^{Q}\alpha_j'p^{j-1}\leq \frac{D}{p-1}\Bigr\}\\
\intertext{and\vspace{-.2em}}
R(D)\,:
={}&\,\max\Bigl\{\textstyle0\leq R\leq\alpha_{Q(D)+1}'\bigm{|}
D_1+\dotsb+D_{\alpha_1'+\dotsb+\alpha_{Q(D)}'+R}\leq \frac{D}{p-1}\Bigr\}\\
={}&\,\max\Bigl\{\textstyle0\leq R\leq\alpha_{Q(D)+1}'\bigm{|}
\sum_{j=1}^{Q(D)}\alpha_j'p^{j-1}+R\,p^{Q(D)}\leq \frac{D}{p-1}\Bigr\}.
\end{align*}
Here, we regard sums of the form $\sum_{j=1}^{0}$ as zero and set
$\alpha_{\alpha_1+1}':=0$ (i.e., $\alpha_{Q(D)+1}':=0$ whenever $Q(D)=\alpha_1$).
Within a graphic representation as in Figure \ref{fig.w} or \ref{fig.w2}, the parameter
$\mu_i(D)$ is the number of red dots in the $i^{\text{th}}$ row, $Q(D)$ is the number of
columns that are completely red, and $R(D)$ is the number of red dots in the next
column, if a next column exists. If there actually is a next column, after the last
completely red one, this column is not completely red, i.e., $R(D)<\alpha_{Q(D)+1}'$ if
(and only if) $Q(D)<\alpha_1$\,.\smallskip
\end{remark}

\begin{remark} \label{rem.V2}\mbox{}
Using the notations in Remark \ref{rem.V}, the minimum value
$\mathcal{V}_p(\underline{\alpha},D)$ can also be expressed in the following forms:
\begin{align*}
\mathcal{V}_p(\underline{\alpha},D)
\,={}&\,\,\sum_{i=1}^{\smash{N}}\bigl(\alpha_i-\mu_i(D)\bigr)\\
\,={}&\,\sum_{i=1}^{\alpha_{Q(D)\!}'}\!\alpha_i\,\,-\alpha_{Q(D)}'Q(D)-R(D)\\
\,={}&\!\sum_{i=1}^{\alpha_{Q(D)+1\!\!\!}'}\!\!\!\alpha_i\,\,-\alpha_{Q(D)+1}'Q(D)-R(D)\\
\,={}&\!\sum_{j=Q(D)+1\!\!\!\!\!\!\!\!}^{\alpha_1}\!\!\!\alpha_j'\,\,-R(D)
\end{align*}
This follows from the formula
$\mathcal{V}_p(\underline{\alpha},D)=\alpha-t(D)=\alpha-\lambda(\underline{\mu}(D))$
at the end of the proof of Theorem \ref{Thm.BC}, and some simple rearrangements. We
illustrated these rearrangements inside Ferrers diagrams in Figure \ref{fig.w2}, where
the positive terms of our formulas are highlighted in green, while the negative terms are
framed in red and blue.
\end{remark}

\noindent As in Corollary \ref{MAINCOR}, our formula for
$\mathcal{V}_p(\underline{\alpha},D)$ simplifies if $\alpha_1=\dotsb=\alpha_N$. With
the parameters
$$Q\,:=\,\bigl\lfloor\log_p(D/N+1)\bigr\rfloor\quad\ \text{and}\ \quad
 R\,:=\,\biggl\lfloor\frac{D-N(p^Q-1)}{(p-1)p^Q}\biggr\rfloor,$$
and with $\overbar{h}:=\max(h,0)$ for real numbers $h$, we obtain the following
corollary:

\begin{cor}\label{Cor.a=1}
If $\alpha_1=\dotsb=\alpha_N$ then
$$\mathcal{V}_p(\underline{\alpha},D)
\,=\, \overbar{N(\alpha_1-Q)-R}.
\medskip$$
If $\alpha_1=\dotsb=\alpha_N=1$ then
$$\mathcal{V}_p(\underline{\alpha},D)
\,=\,\overbar{N-\Bigl\lfloor\frac{D}{p-1}\Bigr\rfloor}.$$
\end{cor}

\begin{proof}
Suppose $\alpha_1=\dotsb=\alpha_N$\,, i.e.,
$\alpha_1'=\dotsb=\alpha_{\alpha_1}'=N$.\smallskip

\noindent\textbf{Case 1, $D<N(p^{\alpha_1}\!-1)$\,:} In this case, it follows as in the
proof of Corollary \ref{MAINCOR} in Section \ref{sec.proof2} that $Q(D)<\alpha_1$
and $R(D)<\alpha_{Q(D)+1}'$. Hence,
$$Q(D)\,=\,\bigl\lfloor\log_p(D/N+1)\bigr\rfloor\,=:\,Q\quad\text{and}\quad
 R(D)\,=\,\biggl\lfloor\frac{D-N(p^Q-1)}{(p-1)p^Q}\biggr\rfloor\,=:\,R\,.$$
So, using the last formula of Remark \ref{rem.V2} (to variate the approach in Section
\ref{sec.proof2}),
$$\mathcal{V}_p(\underline{\alpha},D)
\,=\sum_{j=Q+1\!\!}^{\alpha_1}\!\!\alpha_j'\,\,-R \,=\sum_{j=Q+1\!\!}^{\alpha_1}\!\!N\,\,-R
 \,=\,N(\alpha_1-Q)-R \,=\,\overbar{N(\alpha_1-Q)-R}\,,$$
where the last equality follows from
$N(\alpha_1-Q)-R=\mathcal{V}_p(\underline{\alpha},D)\geq0$\,.\smallskip

In the subcase $\alpha_1=\dotsb=\alpha_N=1$, this further simplifies to
$$\mathcal{V}_p(\underline{\alpha},D)
\,=\,\overbar{N(1-0)-\Bigl\lfloor\frac{D-N(p^0-1)}{(p-1)p^0}\Bigr\rfloor}
\,=\,\overbar{N-\Bigl\lfloor\frac{D}{p-1}\Bigr\rfloor}.$$

\noindent\textbf{Case 2, $D\geq N(p^{\alpha_1}\!-1)$\,:} In this case, $Q\geq
Q(D)=\alpha_1$ and thus $\overbar{N(\alpha_1-Q)-R}=0$\,. By Corollary \ref{Cor.eq},
this is the correct value for $\mathcal{V}_p(\underline{\alpha},D)$ if $D\geq
N(p^{\alpha_1}\!-1)$.\smallskip

The formula for the subcase $\alpha_1=\dotsb=\alpha_N=1$ also gives the correct
value $0$\,.
\end{proof}

\section{Minimization of \protect{$\mathcal N$}}\label{sec.min}

\noindent In this section we determine the minimum
$\min_{\underline{n}\in[\underline{\hat n}(\beta)]}\mathcal N(\underline{n})$ of the
function
$$\mathcal N:[\underline{\hat n}(\beta)]\longrightarrow\N\ ,\quad
\underline{n}\longmapsto\mathcal N(\underline{n}) \,:=\,
\sum_{j=1}^r\overbar{\biggl\lceil\frac{n_j-(p^{\beta_j}\!-1)}{p^{\beta_j-1}(p-1)}\biggr\rceil}
 +\, \mathcal{V}_p\Bigl(\underline{\alpha},\,{\textstyle\sum_{j=1}^r}d_jn_j\Bigr),$$
where (by Theorem \ref{Thm.BC})
$$\mathcal{V}_p\bigl(\underline{\alpha},{\textstyle\sum_{j=1}^r}d_jn_j\bigr)\,=\,\alpha-
 \max\biggl\{0\leq t\leq \alpha\Bigm|
 D_1+\dotsb+D_{t}\leq\tfrac{{\textstyle\sum_{j=1}^r}d_jn_j}{p-1}\biggr\}\,,$$
and where the numbers $\beta,r,\beta_1,\dotsc,\beta_r,d_1,\dotsc,d_r
,N,\alpha_1,\dotsc,\alpha_N\in \Z^+$ with
$$d_1p^{\beta_1}\geq d_2p^{\beta_2}\geq\dotsb\geq d_rp^{\beta_r}
 \quad\text{and}\quad\alpha_1 \geq \alpha_2 \geq \dotsb \geq\alpha_N$$ are
fixed given (and $\beta$ is large enough). Also recall that
$\alpha:=\alpha_1+\alpha_2+\dotsb+\alpha_N$\,, that the numbers
$D_1,D_2,\dotsc,D_{\alpha}$ are defined by
\[\bigl(D_1,D_2,\dotsc,D_{\alpha}\bigr)
\,:=\,\bigl(\,\underbrace{1,1,\dotsc,1}_{\alpha_1'\ \text{times}},
\,\underbrace{p,p,\dotsc,p}_{\alpha_2'\ \text{times}},
\,\dotsc,\,\underbrace{p^{\alpha_1-1},p^{\alpha_1-1},\dotsc,p^{\alpha_1-1}}%
 _{\alpha_{\alpha_1}'\ \text{times}}\,\bigr),\]
that the components of $\hat{\underline{n}}(\beta)=(\hat{n}_j(\beta))_{j=1}^r$ are
given by
$$\hat n_j(\beta) \,\coloneqq\, (p^{\beta_j}\!-1)+(\beta-1)p^{\beta_j-1}(p-1)\,,$$
and that, for every $\underline{\hat n}=(\hat{n}_j)_{j=1}^r\in\N^r$\!,
$$[\underline{\hat n}]:=[\hat n_1]\times[\hat n_2]\times\dotsb\times[\hat n_r]
 \quad\ \text{with}\quad\ [\hat n_j]:=\{0,1,\dotsc,\hat n_j\}\,.$$
\subsection{A Preparatory Lemma}

\noindent It turns out that the minimization of $\mathcal N(\underline{n})$ leads to
another optimization problem that can be stated and solved in more general terms as
follows:

\begin{lemma}
\label{Lem.Smin} Assume $D\in\N$, and let
$\alpha,\Lambda_1,\Lambda_2,\dots,\Lambda_\alpha,V_1,V_2,\dotsc\in\Z^+$\!.
Suppose that $(\Lambda_t)_{t=1}^{\alpha}$ is monotone increasing, that
$(V_t)_{t\in\Z^+}$ is monotone decreasing, and that $\Lambda_1\leq V_1$. Also
assume that $V_{t}=V_{1}$ for all $1\leq t\leq s_0$\,, where
$$s_0 \,:=\, \overbar{\bigl\lceil(\Lambda_1+\dotsb+\Lambda_{t_0}-D)/V_1\bigr\rceil}
 \,\ \text{with}\,\  t_0\,:=\,\max\{1\leq t\leq\alpha\mid \Lambda_t\leq V_1\}\,.\smallskip$$

\noindent Then the function $\mathbf{S}:\N\longrightarrow\Z$ given by
$$\mathbf{S}(s)\,:=\,s-\max\bigl\{0\leq t\leq\alpha\mid
\Lambda_1+\dotsb+\Lambda_{t} \leq V_1+V_2+\dotsb+V_s + D\bigr\}
$$
has a minimum at the point $s_0$\,, and
$$\mathbf{S}(s_0)
\,=\,\begin{cases}
s_0 - t_0 &\text{if $s_0>0$,}\\
- \max\bigl\{0\leq t\leq\alpha\mid
\Lambda_1+\dotsb+\Lambda_{t} \leq D\bigr\} &\text{if $s_0=0$.}
\end{cases}$$
\end{lemma}

\begin{proof}
By definition, $s_0$ is the smallest element of $\N$ with
 $$(\Lambda_1+\dotsb+\Lambda_{t_0}-D)/V_1\leq s_0\,,$$
i.e.,  with
\begin{equation}
\label{eq.s0}\Lambda_1+\dotsb+\Lambda_{t_0}\,\leq\,\,s_0V_1+D\,.
\end{equation}
We calculate $\mathbf{S}(s_0)$, $\mathbf{S}(s_0-s)$ and $\mathbf{S}(s_0+s)$, for all
permissible $s\in\Z^+\!$, to show that $\mathbf{S}(s_0)$ is a minimum of $\mathbf{S}$.
For this purpose it is convenient to extend the sequence $(\Lambda_t)_{t=1}^{\alpha}$
to an infinite sequence by setting
$\Lambda_{\alpha+1},\Lambda_{\alpha+2},\dotsc:=\infty$. With that extension
$t_0=\max\{t\in\Z^+\!\mid \Lambda_t\leq V_1\}$.
\smallskip

\noindent\textbf{Case 1, $s_0>0$\,:} In this case, by \eqref{eq.s0},
\begin{align}
\label{D3}\Lambda_1+\dotsb+\Lambda_{t_0}&\,\leq\,V_1+\dotsb+V_{s_0}+D\\
\intertext{but, by the minimality of $s_0$ in \eqref{eq.s0}, also}
\label{D4}\Lambda_1+\dotsb+\Lambda_{t_0}
 &\,>\,V_1+\dotsb+V_{s_0-1}+D\,.
 \end{align}
In the last inequality, if $s_0\geq2$, each summand $V_j$ on the right is at least as
large as each of the summands $\Lambda_i$ on the left, because
$\Lambda_1\leq\dotsb\leq\Lambda_{t_0}\leq V_1=\dotsb=V_{s_0-1}$. Therefore, we
can remove an equal number of those summands on both sides without destroying the
inequality. Also, the bigger left sum must contain more of the smaller
\(\Lambda\)-summands than the smaller right sum contains of the bigger
\(V\)-summands, because $D\geq0$\,. In particular, for each $0< s\leq s_0$\,,
\begin{equation}
\label{D4+}\Lambda_1+\dotsb+\Lambda_{t_0-s+1}\,>\,V_1+\dotsb+V_{s_0-s}+D\,.
\end{equation}
But, also $\,\dotsb\geq\Lambda_{t_0+2}\geq\Lambda_{t_0+1}>V_1=
 V_{s_0}\geq V_{s_0+1} \geq\dotsb$. So, we can also add an equal number of
subsequent summands on both sides of \eqref{D4}. For each $s\in\N$,
\begin{equation}
\label{D6}\Lambda_1+\dotsb+\Lambda_{t_0+s+1}\,>\,V_1+\dotsb+V_{s_0+s}+D\,.
\end{equation}
Based on these inequalities, we can now calculate $\mathbf{S}(s_0)$,
$\mathbf{S}(s_0-s)$ and $\mathbf{S}(s_0+s)$. It follows from \eqref{D3} and
\eqref{D6} with $s=0$ that
$$\label{S}\mathbf{S}(s_0)\,=\,s_0 - t_0\,.$$
It follows from \eqref{D4+} that, for each $0< s\leq s_0$\,,
$$\mathbf{S}(s_0-s)\,\geq\,s_0 - s - (t_0-s)\,=\,s_0 - t_0\,=\,\mathbf{S}(s_0)\,.$$
And, it follows from \eqref{D6} that, for each $s\in\N$,
$$\mathbf{S}(s_0+s)\,\geq\,s_0 + s - (t_0+s)\,=\,\mathbf{S}(s_0)\,.$$
We see that $\mathbf{S}$ attains a minimum at $s_0$ and
$\mathbf{S}(s_0)=s_0-t_0$\,.\smallskip

\noindent\textbf{Case 2, $s_0=0$\,:} In this case, if we set
$$t(D):=\max\bigl\{0\leq t\leq\alpha\mid \Lambda_1+\dotsb+\Lambda_{t}\leq D\bigr\}
 =\max\bigl\{t\in\N\mid \Lambda_1+\dotsb+\Lambda_{t}\leq D\bigr\}\,,$$
by the maximality of $t(D)$,
\begin{equation}
\label{S0}\Lambda_1+\dotsb+\Lambda_{t(D)+1}\,>\,D\,.
\end{equation}
Moreover,
$\,\dotsb\geq\Lambda_{t(D)+2}\geq\Lambda_{t(D)+1}\geq\Lambda_{t_0+1}>V_1\geq
V_2\geq\dotsb$ since $t(D)\geq t_0$, by \eqref{eq.s0}. Hence, we can add
summands to \eqref{S0}, in the same way as we did it to get \eqref{D6} from
\eqref{D4}. For each $s\in\N$,
$$\Lambda_1+\dotsb+\Lambda_{t(D)+s+1}\,>\,V_1+\dotsb+V_{s}+D\,,$$
and thus
$$\mathbf{S}(s)\,\geq\,s - (t(D)+s)\,=\,0 - t(D)\,=\,\mathbf{S}(0)\,.$$
So, $\mathbf{S}$ attains a minimum at $0$ and $\mathbf{S}(0)= -\max\bigl\{0\leq
t\leq\alpha\mid \Lambda_1+\dotsb+\Lambda_{t}\leq D\bigr\}$.
\end{proof}

\subsection{The Minimum Value of \protect{$\mathcal N$} over
\protect{$[\protect\underline{\hat n}(\beta)]$}}

\noindent We are ready to determine the minimum value of $\mathcal N(\underline{n})$
when $\underline{n}$ is ranging over $[\underline{\hat n}(\beta)]$:

\begin{lemma}
\label{Lem.S} In the settings described at the beginning of Section \ref{sec.min}, with
the derived values
 $\breve\alpha_1,\dotsc,\breve\alpha_N,\breve\alpha,\breve{\mathcal{A}},\mathcal B$
as in Theorem \ref{MAINTHM}, and for every integer
$$\beta\,>\,s_0\,:=\,\overbar{\biggl\lceil\frac{\breve{\mathcal{A}}-\mathcal B}
{d_1p^{{\beta_1-1}}}\biggr\rceil}\,,$$
we have
$$\min_{\underline{n}\in[\underline{\hat n}(\beta)]}\mathcal N(\underline{n})
\,=\,\begin{cases}
s_0+\alpha-\breve\alpha
&\text{if $\breve{\mathcal{A}}>\mathcal B$,}\\
\alpha-\max\bigl\{1\leq t\leq\alpha\mid D_1+\dotsb+D_{t}\leq \mathcal B\bigr\}
&\text{if $\breve{\mathcal{A}}\leq\mathcal B$.}
\end{cases}$$
\end{lemma}

\begin{proof} We shrink the domain $[\underline{\hat n}(\beta)]$ of the variable
$\underline{n}$ till we reach a single point where the minimum is attained and can be
calculated. We proceed in four steps.\smallskip

\noindent\textbf{Step 1:} If $n_1\leq p^{\beta_1}\!-1$ then
$\overbar{\Bigl\lceil\frac{n_1-(p^{\beta_1}\!-1)}{(p-1)p^{\beta_1-1}}\Bigr\rceil}=0$\,. So,
as $\mathcal V_p(\underline{\alpha},\bullet)$ is monotone decreasing,
$$n_1\leq p^{\beta_1}\!-1
\ \ \Longrightarrow\ \ \mathcal N(n_1,n_2,\dotsc,n_r)
 \geq \mathcal N(p^{\beta_1}\!-1,n_2,\dotsc,n_r).
$$
This shows that, in order to find a minimum, we may replace values of $n_1$ below
$p^{\beta_1}\!-1$ with $p^{\beta_1}\!-1\in[\hat n_1(\beta)]$. More generally, for each
$1\le j\leq r$, we may assume $n_j\geq p^{\beta_j}\!-1$. In other words, we may write
each $n_j$ as $u_j+p^{\beta_j}\!-1$ with $u_j\geq0$, which leads to the simplifications
$$\overbar{\biggl\lceil\frac{n_j-(p^{\beta_j}\!-1)}{(p-1)p^{\beta_j-1}}\biggr\rceil}
\,=\, \overbar{\biggl\lceil\frac{u_j}{(p-1)p^{\beta_j-1}}\biggr\rceil}
\,=\, \biggl\lceil\frac{u_j}{(p-1)p^{\beta_j-1}}\biggr\rceil
$$
and
$$
\mathcal{V}_p\Bigl(\underline{\alpha}\,,\,{\textstyle\sum_{j=1}^r}d_jn_j\Bigr)
\,=\,\mathcal{V}_p\Bigl(\underline{\alpha}\,,\,{\textstyle\sum_{j=1}^r}d_ju_j
+ (p-1)\mathcal B\Bigr).
$$
So, with
$$\mathcal U(\underline{u}) \,:=\,
\sum_{j=1}^r\biggl\lceil\frac{u_j}{(p-1)p^{\beta_j-1}}\biggr\rceil + \mathcal{V}_p
 \Bigl(\underline{\alpha}\,,\,{\textstyle\sum_{j=1}^r}d_ju_j + (p-1)\mathcal B\Bigr)$$
we have $\mathcal N(\underline{n})=\mathcal U(\underline{u})$, and thus
$$\min_{\underline{n}\in[\underline{\hat n}(\beta)]}\mathcal N(\underline{n}) \,=\,
 \min_{\underline{u}\in[\underline{\hat u}(\beta)]}\mathcal U(\underline{u})\,.$$
with updated ranges
$$\hat u_j(\beta)\,:=\,\hat n_j(\beta)-(p^{\beta_j}\!-1)\,=\,(\beta-1)p^{\beta_j-1}(p-1)\,.$$
\smallskip

\noindent\textbf{Step 2:} To find a minimum of $\mathcal U$ over $[\underline{\hat
u}(\beta)]=\prod_{j=1}^r[\hat u_j(\beta)]$, we can replace the domain $[\hat
u_j(\beta)]=\{0,1,\dotsc,\hat u_j(\beta)\}$ of each $u_j$ with the smaller domain
\begin{align*}
 [\hat u_j(\beta)]\cap p^{\beta_j -1}(p-1)\Z
 &\,=\,\{0, p^{\beta_j-1}(p-1),\dotsc, (\beta-1)p^{\beta_j-1}(p-1)\}\\
 &\,=\,p^{\beta_j-1}(p-1)[\beta-1].
\end{align*}
Indeed, if the $j^{\text{th}}$ argument $u_j\in[\hat u_j]$ of $\mathcal U(u_1,\dotsc,u_r)$
is replaced with the first multiple of $p^{\beta_j-1}(p-1)$ above or equal to $u_j$ (which
still lies inside $[\hat u_j(\beta)] = [(\beta-1)p^{\beta_j-1}(p-1)]$), then the summand
$\Bigl\lceil\frac{u_j}{p^{\beta_j-1}(p-1)}\Bigr\rceil$ of $\mathcal U(\underline{u})$ stays
the same and $\mathcal U(\underline{u})$ certainly does not increase. The minimum of
$\mathcal U$ is already attained at a point $\underline{u}$ of the smaller domain
$\prod_{j=1}^r\bigl(p^{\beta_j-1} (p-1)[\beta-1]\bigr)\subseteq[\underline{\hat u}(\beta)]$.
Hence, with
\begin{align*}
\mathcal T(t_1,\dotsc,t_r)
\ :={}&\ \ \mathcal U\bigl(p^{\beta_1-1}(p-1)t_1,\dotsc,p^{\beta_r-1}(p-1)t_r\bigr)\\
={}&\ \ t_1  + \dotsb + t_r
+ \mathcal{V}_p\Bigl(\underline{\alpha}\,,\,{\textstyle\sum_{j=1}^r}d_jp^{\beta_j-1}(p-1)t_j
+ (p-1)\mathcal B\Bigr)
\end{align*}
we have
$$\min_{\underline{n}\in[\underline{\hat n}]}\mathcal N(\underline{n}) \,=\,
\min_{\underline{u}\in[\underline{\hat u}]}\mathcal U(\underline{u}) \,=\,
\min_{\underline{t}\in[\beta-1]^r}\mathcal T(\underline{t}).\smallskip$$

\noindent\textbf{Step 3:} In our search for the minimum value that $\mathcal
T(t_1,\dotsc,t_r)$ may take, we can now modify any two arguments $t_i$ and $t_j$ with
$i<j$ by replacing $t_j$ with $t_j-1$ and $t_i$ with $t_i+1$. If we view the term
$d_jp^{\beta_j-1}(p-1)t_j$ as sum of $t_j$ equal summands $d_jp^{\beta_j-1}(p-1)$, this
step changes one of the $t_j$ summands $d_jp^{\beta_j-1}(p-1)$ inside the argument of
$\mathcal V_p(\underline{\alpha},\bullet)$ into one additional summand
$d_ip^{\beta_i-1}(p-1)$, of which we then have $t_i+1$. Since
\[d_1p^{\beta_1-1}\geq d_2p^{\beta_2-1}\geq\dotsb\geq d_rp^{\beta_r-1} \] and
$\mathcal V_p(\underline{\alpha},\bullet)$ is monotone decreasing, we have
$$\mathcal T(\dotsc,t_i+1,\dotsc,t_j-1,\dotsc)
 \,\leq\,\mathcal T(\dotsc,t_i,\dotsc,t_j,\dotsc).$$
The only restriction to such modifications is that all argument $t_j$ must stay within their
domains $[\beta-1]$. They cannot increase above $\beta-1$ or go below $0$. Through
repeated applications of our modification, we can empty some $t_j$ and fill others. This
shows that the minimum is attained at a point of the form
$$(t_1,t_2,\dotsc,t_r)\,=\,(\beta-1,\beta-1,\dotsc,\beta-1,x,0,0,\dotsc,0).$$
At such points, we have
$$\mathcal T(t_1,t_2,\dotsc, t_r)
\,=\,s + \mathcal{V}_p\Bigl(\underline{\alpha}\,,\,(p-1)(V_1+V_2+\dotsb+V_s
 + \mathcal B)\Bigr),$$
where $s=t_1+t_2+\dotsb+t_r=\beta-1+\beta-1+\dotsb+\beta-1+x\leq r(\beta-1)$, and
where
$$
\bigl(V_1,V_2,\dotsc,V_{r(\beta-1)}\bigr)\,:=\,
\bigl(\,\underbrace{d_1p^{\beta_1-1}\!,
\dotsc,d_1p^{\beta_1-1}}_{\beta-1\ \text{times}\!},
\,\dotsc, \,\underbrace{d_rp^{\beta_r-1}\!,\dotsc,d_rp^{\beta_r-1}}_{\beta-1\
 \text{times}}\,\bigr).$$
Hence, with the function
$$\mathcal S:[r(\beta-1)]\ra\N\,,\quad\mathcal S(s)\,:=\,
s + \mathcal{V}_p\Bigl(\underline{\alpha}\,,\,(p-1)(V_1+V_2+\dotsb+V_s
 + \mathcal B)\Bigr),$$
we have
$$\min_{\underline{n}\in[\underline{\hat n}(\beta)]}\mathcal N(\underline{n}) \,\,=
\min_{\underline{t}\in[\beta-1]^r}\mathcal T(\underline{t}) \,\,=
\min_{s\in[r(\beta-1)]}\mathcal S(s).\smallskip$$

\noindent\textbf{Step 4:} To find the minimum of $\mathcal S$, we use Lemma
\ref{Lem.Smin} with $D:=\mathcal B$, $\alpha:=\alpha_1+\dotsb+\alpha_N$, and
$\Lambda_t:=D_t$ for all $1\leq t\leq\alpha$. We also use the values $V_t$ as defined
above for all $1\leq t\leq r(\beta-1)$, and set $V_t:=V_{r(\beta-1)}$ for all $t>r(\beta-1)$.
With the infinite sequence $(V_t)_{t\in\Z^+}$ the domain of $\mathcal S$ can be
extended to $\N$ (with the hope not to alter its minimum in doing so), as the expression
$$\mathcal{V}_p\Bigl(\underline{\alpha}\,,\,(p-1)(V_1+\dotsb+V_s + \mathcal B)\Bigr)
 \,=\,\alpha-\max\bigl\{0\leq t\leq \alpha\mid\Lambda_1+\dotsb+\Lambda_{t}\leq
 V_1+\dotsb+V_s + \mathcal B\bigr\}$$
makes sense for all $s\in\N$. The extended function $\mathcal S:\N\ra\N$ is then
almost the same as the function $\mathbf{S}:\N\ra\Z$ in Lemma \ref{Lem.Smin}. For all
$s\in\N$,
$$\mathcal S(s)\,:=\,\mathbf{S}(s) + \alpha\,.$$
We also have $\Lambda_1\leq V_1$ as required in Lemma \ref{Lem.Smin}.
Moreover, as in our situation the sequence $(\Lambda_t)$ contains repetitions of
lengths $\alpha_1',\alpha_2',\dotsc,\alpha'_{\alpha_1}$, the parameter
$$t_0\,:=\,\max\{1\leq t\leq\alpha\mid\Lambda_t\leq V_1\}$$
in Lemma \ref{Lem.Smin} can be written as
$$t_0\,=\,\alpha_1'+\dotsb+\alpha_{i_0}'\quad\text{with}\quad
 i_0\,:=\,\max\{1\leq i\leq\alpha_1\mid p^{i-1}\leq d_1p^{\beta_1-1}\}\,.$$
Here, the inequality $p^{i-1}\leq d_1p^{\beta_1-1}$ can be written as
$i\leq\beta_1+\log_p\bigl(d_1\bigr)$, and the biggest integer $i$ with this property is
$\beta_1+\lfloor\log_p\bigl(d_1\bigr)\rfloor$. But, in the definition of $i_0$ we also have
the requirement $i\leq\alpha_1$, so that
$$i_0\,=\,\max\{1\leq i\leq\alpha_1\mid p^{i-1}\leq d_1p^{\beta_1-1}\}
 \,=\,\min\bigl\{\alpha_1,\beta_1+\lfloor\log_p\bigl(d_1\bigr)\rfloor\bigr\}
 \,=\,\breve\alpha_1\,.$$
With that and the second equation of Example \ref{Ex.con2} we get
$$t_0\,=\,\alpha_1'+\dotsb+\alpha_{i_0}'
\,=\,\alpha_1'+\dotsb+\alpha_{\breve\alpha_1}'
\,=\,\breve\alpha_1+\dotsb+\breve\alpha_N
 \,=\,\breve\alpha\,.$$
With the third equation of Example \ref{Ex.con2}, we further see that
$$\Lambda_1+\dotsb+\Lambda_{t_0}
 \,=\,\Lambda_1+\dotsb+\Lambda_{\breve\alpha}
 \,=\,\alpha_1'p^0+\dotsb+\alpha_{\breve\alpha_1}'p^{\breve\alpha_1-1}
 \,=\,\sum_{i=1}^N\frac{p^{\breve\alpha_i}\!-1}{p-1}
 \,=:\,\breve{\mathcal{A}}\,.$$
In particular, the definition of $s_0$ in Lemma \ref{Lem.Smin} coincides with the
current one:
$$s_0
 \,=\,\overbar{\bigl\lceil(\Lambda_1+\dotsb+\Lambda_{t_0}-\mathcal B)/V_1\bigr\rceil}
 \,=\,\overbar{\biggl\lceil\frac{\breve{\mathcal{A}}
 -\mathcal B}{d_1p^{{\beta_1-1}}}\biggr\rceil}
 \,.$$
As assumed, this number is smaller than $\beta$, i.e., $s_0\leq\beta-1$. This shows
that $V_{t}=V_{1}$ for all $1\leq t\leq s_0$, as required in Lemma \ref{Lem.Smin}. But,
$s_0\leq\beta-1$ also shows that the minimum point $s_0$ of $\mathbf{S}$ lies inside
$[r(\beta-1)]$. Hence, the minimum point $s_0$ of $\mathbf{S}$ is also a minimum point
of $\mathbf{S}|_{[r(\beta-1)]}$ and of $\mathcal S|_{[r(\beta-1)]}$. Thus, Lemma
\ref{Lem.Smin} yields
\begin{align*}
\min_{\underline{n}\in[\underline{\hat n}(\beta)]}\mathcal N(\underline{n})
 &\,=\,\min_{s\in[r(\beta-1)]}\mathcal S(s)\\
 &\,=\,\,\,\mathbf{S}(s_0)\,+\,\alpha\\
 &\,=\,\begin{cases}
s_0+\alpha-\breve\alpha
&\text{if $\breve{\mathcal{A}}>\mathcal B$,}\\
\alpha-\max\bigl\{1\leq t\leq\alpha\mid D_1+\dotsb+D_{t}\leq \mathcal B\bigr\}
&\text{if $\breve{\mathcal{A}}\leq\mathcal B$.}
\end{cases}
\end{align*}
\end{proof}

\section{About Conjugate Partitions}\label{sec.CS}

\noindent
\newcommand{\al}a%
In the previous sections, we repeatedly used results about conjugate sequences. In this
section, we discus and prove those results in the form of lemmas and examples. Given
a sequence $(\al_1,\al_2,\dotsc,\al_N)$ of integers with $N>0$ and
$\al_1\geq\al_2\geq\dotsb\geq\al_N>0$, the conjugate numbers
$\al_{1}',\al_{2}',\dotsc,\al_{\al_1}'$ are defined by
$$\al_j'\,:=\,\#\bigl\{1\leq i\leq N\mid\al_i\geq j\bigr\}\,.$$
The sequence $\underline{a}:=(\al_i)\in(\Z^+)^{N}$ is a partition of the number
$\al:=\al_1+\al_2+\dotsb+\al_N$, but the finite monotone decreasing sequence
$\underline{a}':=(\al_j')\in(\Z^+)^{\al_1}$ also partitions $\al$ (i.e.,
$\al:=\al_1'+\al_2'+\dotsb+\al_{\al_1}'$), as we will see. It is called the \textbf{conjugate
partition}, and it is a dual partition of $\al$, in the sense that the conjugate of the
conjugate is the original sequence. This is easy to see if we represent each $a_i$ by a
row of $a_i$ dots, in a so called \textbf{Ferrers diagram}. The conjugate partition is then
obtained by reflecting the corresponding Ferrers diagram about the main diagonal, like
transposing a matrix:

\begin{figure}[h!]
\begin{center}
\includegraphics[width=0.65\textwidth]{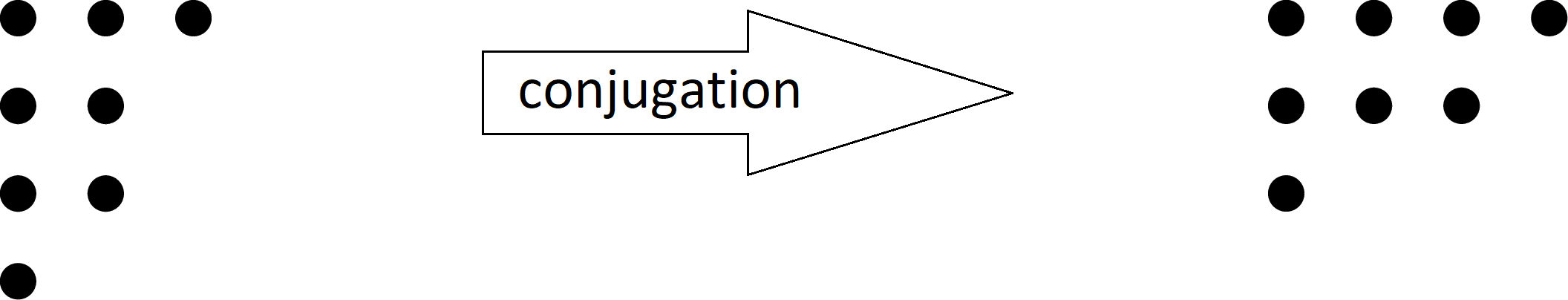}
\end{center}
\caption{The conjugate of $(3,2,2,1)$ is $(4,3,1)$.}
\end{figure}

\subsection{Two General Conjugation Lemmas}

\noindent The following lemma is formulated in a way that is helpful in our
calculations.

\begin{lemma}\label{Lem.con}
Let $(a_i)\in(\Z^+)^N$ be monotone decreasing, and $1\leq m\leq a_1$. We have the
following identity in $\Z[x]$:
$$
a_{m}'x^{m}+a_{m+1}'x^{m+1}+\dotsb+a_{a_1}'x^{a_1}
 \,=\,\sum_{i=1}^{a_{m}'}\bigl(x^{m}+x^{m+1}+\dotsb+x^{a_i}\bigr)\,.
$$
\end{lemma}

\begin{proof}
Both polynomials have degree at most $a_1$, and there are no monomials of degree
less than $m$. For each $m\leq j\leq a_1$, however, the coefficient of $x^j$ in the
standard expansion of the right polynomial is
$$\#\bigl\{1\leq i\leq a_{m}'\mid a_i\geq j\bigr\}
 \,=\,\#\bigl\{1\leq i\leq N\mid a_i\geq j\bigr\}\,=\,a_j'\,,$$
because
$$a_{i}\geq j\quad\Longrightarrow\quad a_{i}\geq m\quad\Longrightarrow\quad
 a_{1},\dotsc,a_{i}\geq m\quad\Longrightarrow\quad i\leq a_{m}'\,,$$
i.e., it is the same as the coefficient of $x^j$ in the left polynomial.
\end{proof}

\noindent The following lemma is clear if we imagine taking the minimum as intersecting
two Ferrers diagrams, because ``intersecting'' and ``reflecting'' commute.

\begin{lemma}\label{Lem.con2}
If the two sequences $(a_i),(b_i)\in(\Z^+)^N$ are monotone decreasing, then the
sequence $(c_i):=\bigl(\min(a_i,b_i)\bigr)\in(\Z^+)^N$ is also monotone decreasing.
The conjugate sequences $(a_j')$, $(b_j')$, and $(c_j')$ have lengths $a_1$, $b_1$,
and $c_1=\min(a_1,b_1)$, respectively; and for all $1\leq j\leq c_1$,
$$c_j'\,=\,\min(a_j',b_j')\,.$$
\end{lemma}

\subsection{Special Cases}

\begin{example}\label{Ex.con1}
If $x=1$ in Lemma \ref{Lem.con}, we obtain, for all $1\leq m\leq a_1$,
$$\sum_{j=m}^{a_1}\! a_j'\,=\,\sum_{i=1}^{a_{m}'}\bigl(a_i-m+1\bigr)\,.$$
If $m=1$ then $a_m'=a_1'=N$, and then $x=1$ leads to
$$a_1'+\dotsb+a_{a_1}'\,=\,a_1+\dotsb+a_{N}\,.$$
If $m=1$ and $x=p$, then Lemma \ref{Lem.con} yields
$$a_1'p^0+\dotsb+a_{a_1}'p^{a_1-1}
\,=\,\sum_{i=1}^{N}\bigl(p^{0}+p^{1}+\dotsb+p^{a_i-1}\bigr)
 \,=\,\sum_{i=1}^N\frac{p^{a_i}\!-1}{p-1}\,.$$
\end{example}

\begin{example}\label{Ex.con2}
If the sequence $(b_i)$ in Lemma \ref{Lem.con2} is constant equal to a number
$b\in\Z^+$\!, we obtain as the conjugate of the sequence
$(c_i):=\bigl(\min(a_i,b)\bigr)\in(\Z^+)^N$ the sequence
$$(c_1',c_2',\dotsc,c_{c_1}')\,=\,(a_1',a_2',\dotsc,a_{c_1}')\,.$$
So, if we apply Lemma \ref{Lem.con} to $(c_i)$ with $m=1$ and $x=1$, we get
$$
a_1'+\dotsb+a_{c_1}'\,=\,c_1+\dotsb+c_{N}\,.
$$
If instead $m=1$ and $x=p$, we get
$$
a_1'p^0+\dotsb+a_{c_1}'p^{c_1-1}
\,=\,\sum_{i=1}^{N}\bigl(p^{0}+p^{1}+\dotsb+p^{c_i-1}\bigr)
\,=\,\sum_{i=1}^N\frac{p^{c_i}\!-1}{p-1}\,.
$$
\end{example}

\end{document}